\documentclass[12pt]{article}
\usepackage[utf8]{inputenc}
\usepackage{amssymb,amsmath,amsfonts,amsthm,amscd,latexsym,indentfirst,verbatim}
\usepackage[english]{babel}
\usepackage{geometry,color,ulem}
\def\dbl{\lbrace\kern-3pt\lbrace}
\def\dbr{\rbrace\kern-3pt\rbrace}

\def\Imm{\textrm{Im}}

\def\Aut{\textrm{Aut}}

\def\Spec{\mathrm{Spec}}

\begin{document}
\newcommand{\tocsecindent}{\hspace{7mm}}

\sloppy
\hfill{20D15 (MSC2020)}

\begin{center}
{\Large
Rota---Baxter operators on extraspecial groups}

\smallskip

Andrey Savelyev
\end{center}

\begin{abstract}
Theory of Rota---Baxter operators on rings and algebras has been developed
since 1960. In 2021, L. Guo, H. Lang, Y. Sheng have defined the
notion of Rota—Baxter operator on a group. 
We classify all Rota---Baxter operators on the finite Heisenberg group and the extraspecial group of order~$p^3$. We also investigate an analogue of the spectral property for Rota---Baxter operators on algebras.

\end{abstract}

\section{Introduction}
Rota---Baxter operators on algebras, which are an algebraic analogue of integration operators, have been studied since the 1960s. By now, connections have been found between Rota---Baxter operators on algebras and the Yang---Baxter equation, pre- and post-algebras, double Lie algebras, Bernoulli polynomials, etc., see~\cite{GuoMonograph}. The theory has received its further development, including the definition of Rota---Baxter operators on other algebraic structures, for example, on groups and Hopf algebras.

\textbf{Definition 1}~\cite{Guo}.
Let $G$ be a group, a mapping $B\colon G \to G$ is called a Rota---Baxter operator on the group $G$, if for any $g, h \in G$ the equality holds
\begin{equation} \label{RBGroup}
B(g)B(h) = B(g B(g) h B(g)^{-1}).    
\end{equation}
A group $G$ equipped with a Rota---Baxter operator $B$ is called a Rota---Baxter group. 
It was shown in~\cite{Guo} that differentiating a continuous Rota---Baxter operator on a Lie group yields a Rota---Baxter operator of weight 1 on the corresponding Lie algebra.
A~connection between Rota---Baxter groups, braces, and solutions to the quantum Yang–Baxter equation was established in~\cite{BG}.

One of the directions of the study of Rota---Baxter operators is their classification on particular algebraic structures, for example \cite{BGspor, GG, Goncharov, GP, PBG}. In the context of group Rota---Baxter operators, papers \cite{BGspor} and \cite{GG} are of interest, where Rota---Baxter operators on simple sporadic, dihedral, and alternating groups are classified. In this paper, Rota---Baxter operators on the finite Heisenberg group $H_3(\mathbb{Z}_p)$ and the extraspecial group of order $p^3$, $p \geq 3$, are classified.

Prior to this, a series of papers~\cite{Ita,Goncharovplus,Zhu} considered various examples of Rota---Baxter operators on Heisenberg groups $H_3(\mathbb{F})$, 
where $\mathbb{F} = \mathbb{Z}_p$, $\mathbb{R}$, or $\mathbb{C}$.
In the paper~\cite{Rathee}, N.~Rathee showed how 
to define a Rota---Baxter operator on the group $H_3(\mathbb{C})$
from a~Rota---Baxter operator of weight 1 on the Heisenberg Lie algebra $\mathfrak{h}_3(\mathbb{C})$. We note that Rota---Baxter operators of weight 1 on the algebra $\mathfrak{h}_3(\mathbb{C})$
were described in~\cite{Ji}.
At the same time, an approach generalizing the results of N. Rathee was presented in the paper~\cite{Goncharovplus}. In connection with the appearance of a large number of examples of operators in the literature, the problem of obtaining a complete classification of Rota---Baxter operators on a finite Heisenberg group seems extremely interesting. The extraspecial group of order $p^3$, in its turn, is an analogue of a finite Heisenberg group; these are the only two non-isomorphic groups of order 27 in the case $p = 3$.

Therefore, a complete classification of Rota---Baxter operators on the finite Heisenberg group is of significant interest due to many examples appearing in the literature. The extraspecial group of order~$p^3$ is an analogue of the finite Heisenberg group. For~$p = 3$, these two structures constitute the only non-isomorphic groups of order~27.

Returning to the connection between Rota---Baxter operators and braces, we note the 2025 Ph.D. thesis by S. Trappeniers~\cite{Trappeniers}. There, the author described non-abelian braces on the finite Heisenberg group $H_3$ and the extraspecial group $E_3$ of order $p^3$. In addition, extraspecial braces were studied in a recent paper~\cite{BEKP}.

In addition to the classification, this paper proves that Rota---Baxter operators on any extraspecial group preserve the center. 
Then we analyze the obtained results.

In Section~2, we provide the required preliminaries on Rota---Baxter operators on groups and automorphisms of $H_3$ and $E_3$.

In Section~3, we prove that Rota---Baxter operators preserve the center of any extraspecial group. 
Namely, for an arbitrary extraspecial group~$G$ and any Rota---Baxter operator~$B$ on it, $B(Z(G)) \subseteq Z(G)$, where~$Z(G)$ is the center of~$G$.

In Section~4, we apply the results from Section~3 to determine the general form of Rota---Baxter operators on the finite Heisenberg group and the extraspecial group of order~$p^3$. 
Then we completely classify these operators.

In Section~5, we analyze the obtained results. 
In particular, we determine which operators are splitting, investigate an analogue of the spectral property, and study bijective Rota---Baxter operators on the Heisenberg group in this context.

\section{Preliminaries}
\subsection{Definitions and properties}
\textbf{Definition 2}.
$p$-group $G$  is called extraspecial if $|Z(G)| = p$ and
$G/Z(G) \cong \mathbb{Z}_p \oplus \ldots \oplus \mathbb{Z}_p$.

Let $G$ be group, define $[g, h] = g^{-1} h^{-1} g h$.

\textbf{Definition 3}.
A group $H_3 = \langle a, b, c \mid [a, b] = c, \, [a, c] = [b, c] = a^p = b^p = c^p = e \rangle$ is called the finite Heisenberg group of order $p^3$. 

\textbf{Definition 4}.
A group $E_3 = \langle a, b, c \mid [a, b] = c, \, a^p = b^p = c, \, [a, c] = [b, c] = c^p = e \rangle$ is called extraspecial group of order $p^3$.

\textbf{Definition 5}.
Let $G$ and $H$ be groups and $\psi \colon Z(G) \to Z(H)$ be isomorphism. 
Then the central product of groups $G$ and $H$ is the quotient group
$(G \times H) / K$, where 
\begin{equation} \label{center_product}
K = \{ (g, \psi(g)) \in Z(G) \times Z(H)\} \trianglelefteq G \times H.
\end{equation}
We will denote $(G \times H) / K$ as $G \times_\psi H$.

Let us discuss an equivalent definition of extraspecial groups.

\textbf{Remark 1} \cite{Craven}.
Any extraspecial group of order $p^{2n+1}$ with $n \in \mathbb{N}$ is isomorphic to the central product of $n$ copies of $H_3$ or to the central product of $n-1$ copies of $H_3$ and one copy of $E_3$.

Now we will note some general properties for Rota---Baxter groups.

\textbf{Proposition 1}~\cite{Guo}. 
Let $B$ be a Rota---Baxter operator on a group $G$. Then

a) $\Imm(B),\ker(B)\leq  G$,

b) $B(e) = e$.

We define a map $\widetilde{B}\colon G\to G$ for Rota---Baxter operator $B$ and $G$ as follows:
$$
\widetilde{B}(g) = g^{-1}B(g^{-1}), \quad g \in G.
$$ 
It is easy to see that $\widetilde{\widetilde{B}}(g) = B$.

\textbf{Proposition 2}~\cite{Guo}. 
Let $B$ be a Rota---Baxter operator on a group $G$. Then

a) The map $\widetilde{B}$ is a Rota---Baxter operator on $G$.

b) If $|G| < \infty$, then $|G| = |\ker(B)| \cdot |\Imm(B)|$.

c) The map $\varphi^{-1} B \varphi$ is a Rota---Baxter operator on $G$ for any $\varphi \in \Aut(G)$.

A group $G$  is called factorizable if $G = HL$ for some subgroups~$H,L$.
The expression $G = HL$ is called the factorization of the group $G$.
If $H\cap L = \{e\}$, then such a~factorization is called exact.

{\bf Example 1}~\cite{Guo}.
Let $G$ be group and we have exact factorization $G = HL$.
Then a~map $B\colon G\to G$ defined by the formula
$B(hl) = l^{-1}$ is a RB-operator on~$G$.

\textbf{Definition 6}~\cite{Guo}.
An operator defined as shown in Example 1 is called a splitting RB-operator.

\textbf{Proposition 3}~\cite{BG}. 
An operator $B$ on a group $G$ is splitting if and only if $B\widetilde{B}(G) = \{e\}$.

We note some connection between the images and kernels of operators~$B$ and $\widetilde{B}$.

\textbf{Proposition 4}~\cite{GG,Guo}.
Let $B$ be Rota---Baxter operator on a group $G$. Then

a) $\ker(B) \trianglelefteq \Imm(\widetilde{B})$,

b) $G = \Imm(B) \cdot \Imm(\widetilde{B})$,

c) $\Imm(B \widetilde{B}) = \Imm(\widetilde{B} B) = \Imm(B) \cap \Imm(\widetilde{B})$.

For brevity, we will write $\bar{k}\in \mathbb{Z}_{p^2}$ or $\bar{k} \in \mathbb{Z}_p$ simply as~$k$ from this moment.

For an element $k \in \mathbb{Z}_{p^2}$ unique $r, q \in \{ 0, \ldots, p-1 \} \subseteq \mathbb{Z}_{p^2}$ exist such that $k = r + pq$. Define $k \ \mathrm{rem} \ p = r$.

\subsection{Automorphisms of $H_3$ and $E_3$}

In this section, we describe the automorphisms of the groups $H_3$ and $E_3$.
The description of the structures of the groups $\Aut(H_3)$ and $\Aut(E_3)$ can be found in \cite{Winter}. For convenience of further use, we give explicit formulas for the automorphisms of these groups in terms of the introduced copresentations of these groups.

Next, we define $\varphi_\alpha, \varphi_\beta \colon \mathbb{Z}_p \times \mathbb{Z}_p \to \mathbb{Z}_p$, $\varphi_\gamma \colon \mathbb{Z}_p \times \mathbb{Z}_p \times \mathbb{Z}_p \to \mathbb{Z}_p$ for $G = H_3$ and $\varphi_\alpha, \varphi_\beta \colon  \{0,\ldots,p-1\} \times \{0,\ldots,p-1\} \to \mathbb{Z}_{p^2}$, $\varphi_\gamma \colon  \{0,\ldots,p-1\} \times \{0,\ldots,p-1\} \times \{0,\ldots,p-1\} \to \mathbb{Z}_{p^2}$ for $G = E_3$.

Let $0 \leq t, s, r \leq p-1$, $\varphi \in \Aut(G)$, then
\begin{equation} \label{define_automorphism}
\varphi(b^s a^t c^r) = b^{\varphi_\beta(t, s)} a^{\varphi_\alpha(t, s)} c^{\varphi_\gamma(t, s, r)}.
\end{equation}
We have
\begin{equation} \label{phi_alpha_beta}
\varphi_\alpha(t, s) = u_\alpha t + u_\beta s, \quad \varphi_\beta(t, s) = v_\alpha t + v_\beta s, \quad \varphi_\gamma(t, s, r) = \varphi_\gamma(t, s, 0) + \Delta r,
\end{equation}
since $\varphi$ is automorphism, where $u_\alpha, u_\beta, v_\alpha, v_\beta \in \mathbb{Z}_p$, $\Delta \in \mathbb{Z}_p^{*}$, $u_\alpha v_\beta - u_\beta v_\alpha \neq 0$ if $G = H_3$ and $u_\alpha, u_\beta, v_\alpha, v_\beta \in \mathbb{Z}_{p^2}$, $\Delta \in \mathbb{Z}_{p^2}^{*}$, $u_\alpha v_\beta - u_\beta v_\alpha \neq 0 \pmod{p}$ if $G = E_3$.

The following relations hold:
\begin{gather*}
b^{v_\alpha t + v_\beta s} a^{u_\alpha t + u_\beta s} c^{\varphi_\gamma(t, s, 0)} = \varphi( b^s a^t) =\varphi(b^s)  \varphi(a^t) = b^{v_\alpha t + v_\beta s} a^{u_\alpha t + u_\beta s} c^{\varphi_\gamma(t, 0, 0) + \varphi_\gamma(0, s, 0) + u_\beta v_\alpha t s},
\\
b^{v_\alpha + v_\beta} a^{u_\alpha + u_\beta} c^{\varphi_\gamma(1, 1, 0) + \Delta} 
 {=} \varphi(b a c) 
 {=} \varphi(a b) 
 {=} \varphi(a) \varphi(b) 
 {=} b^{v_\alpha + v_\beta} a^{u_\alpha + u_\beta} c^{\varphi_\gamma(1, 0, 0) + \varphi_\gamma(0, 1, 0) + u_\alpha v_\beta}.
\end{gather*}

Therefore, $\varphi_\gamma(t, s, 0) = \varphi_\gamma(t, 0, 0) + \varphi_\gamma(0, s, 0) + u_\beta v_\alpha t s \pmod{p}$ and $\Delta = u_\alpha v_\beta - u_\beta v_\alpha \pmod{p}$.

Let us consider the action on the powers of the element
\begin{gather*}
b^{v_\alpha t} a^{u_\alpha t} c^{\varphi_\gamma(t, 0, 0)} = \varphi(a^t) = \varphi(a)^t = b^{v_\alpha t} a^{u_\alpha t} c^{t \varphi_\gamma(1, 0, 0) + \frac{t(t-1)}{2} u_\alpha v_\alpha}, \\
b^{v_\beta s} a^{u_\beta s} c^{\varphi_\gamma(0, s, 0)} = \varphi(b^s) = \varphi(b)^s = b^{v_\beta s} a^{u_\beta s} c^{s \varphi_\gamma(0, 1, 0) + \frac{s(s-1)}{2} u_\beta v_\beta}.
\end{gather*}

Denote $w_\alpha = \varphi_\gamma(1, 0, 0)$, $w_\beta = \varphi_\gamma(0, 1, 0)$, then $\varphi_\gamma(t, 0, 0) = w_\alpha t + \frac{t(t-1)}{2} u_\alpha v_\alpha$ and $\varphi_\gamma(0, s, 0) = w_\beta s + \frac{s(s-1)}{2} u_\beta v_\beta$.

We obtain the general formula for $\varphi_\gamma$:
\begin{equation} \label{phi_gamma}
\varphi_\gamma(t, s, 0) = w_\alpha t + w_\beta s + \frac{t(t-1)}{2} u_\alpha v_\alpha + \frac{s(s-1)}{2} u_\beta v_\beta + u_\beta v_\alpha t s.
\end{equation}

We now determine what conditions do we need on the coefficients $u_\alpha, u_\beta, v_\alpha, v_\beta, w_\alpha$, $w_\beta$, for the mapping $\varphi$ to be an automorphism.

Note that $\varphi$ is bijective if and only if $\Delta = u_\alpha v_\beta - u_\beta v_\alpha \neq 0 \pmod{p}$. Let us check the homomorphism condition:
\begin{multline*}
b^{v_\alpha(t_1 + t_2 \ \mathrm{rem} \ p) + v_\beta(s_1 + s_2 \ \mathrm{rem} \ p)} a^{u_\alpha(t_1 + t_2 \ \mathrm{rem} \ p) + u_\beta(s_1 + s_2 \ \mathrm{rem} \ p)} c^{C_1} \\
 = \varphi(b^{s_1 + s_2 \ \mathrm{rem} \ p} a^{t_1 + t_2 \ \mathrm{rem} \ p} c^{r_1 + r_2 + t_1 s_2 + R(t_1 + t_2, s_1 + s_2)}) 
 = \varphi(b^{s_1} a^{t_1} c^{r_1}) \varphi(b^{s_2} a^{t_2} c^{r_2}) \\
 = b^{v_\alpha(t_1 + t_2 \ \mathrm{rem} \ p) + v_\beta(s_1 + s_2 \ \mathrm{rem} \ p)} a^{u_\alpha(t_1 + t_2 \ \mathrm{rem} \ p) + u_\beta(s_1 + s_2 \ \mathrm{rem} \ p)} c^{C_2},
\end{multline*}
where the following formulas hold modulo $p$:
\begin{multline*}
C_1 = \varphi_\gamma(t_1 + t_2 \ \mathrm{rem} \ p, s_1 + s_2 \ \mathrm{rem} \ p, 0) + \Delta (r_1 + r_2 + t_1 s_2 + R(t_1 + t_2, s_1 + s_2))\\
= w_\alpha (t_1 + t_2) + w_\beta (s_1 + s_2) + \frac{(t_1 + t_2)(t_1 + t_2 -1)}{2} u_\alpha v_\alpha + \frac{(s_1 + s_2)(s_1 + s_2 - 1)}{2} u_\beta v_\beta \\
+ u_\beta v_\alpha (t_1 + t_2)(s_1 + s_2) + \Delta (r_1 + r_2 + t_1 s_2 + R(t_1 + t_2, s_1 + s_2));
\end{multline*}

\vspace{-0.8cm}

\begin{multline*}
C_2 = w_\alpha (t_1 + t_2) + w_\beta (s_1 + s_2) + \bigg( \frac{t_1(t_1-1)}{2} + \frac{t_2(t_2-1)}{2} \bigg) u_\alpha v_\alpha \\
+ \bigg( \frac{s_1(s_1-1)}{2} + \frac{s_2(s_2-1)}{2} \bigg) u_\beta v_\beta + u_\beta v_\alpha (t_1 s_1 + t_2 s_2) + \Delta (r_1 + r_2) + Q(t_1 + t_2, s_1 + s_2) \\
+ (u_\alpha t_1 + u_\beta s_1)(v_\alpha t_2 + v_\beta s_2);
\end{multline*}
$$ 
R(t, s) = [t/p] + [s/p], \quad 
Q(s, t) = (u_\alpha + v_\alpha)[t/p] + (u_\beta + v_\beta)[s/p] \colon \mathbb{Z}_{p^2} \times \mathbb{Z}_{p^2} \to \mathbb{Z}_{p^2}.
$$

We open the brackets and use the fact $\Delta = u_\alpha v_\beta - u_\beta v_\alpha$, then $C_1 = C_2 \pmod{p}$ if and only if
$\Delta R(t, s) = Q(t, s) \pmod{p}$. 
This condition is fulfilled automatically for  $G = H_3$ and for $G = E_3$ it is equivalent to the equalities
$u_\alpha + v_\alpha = u_\beta + v_\beta = \Delta \pmod{p}$.
Thus, for $G = E_3$, we have the following equalities modulo $p$:
$$
0 \neq \Delta 
 = u_\alpha v_\beta - u_\beta v_\alpha 
 = u_\alpha (\Delta - u_\beta) - u_\beta (\Delta - u_\alpha) 
 = \Delta (u_\alpha - u_\beta),
$$
Hence, $u_\beta = u_\alpha - 1 \pmod{p}$.

Let $\varphi \in \Aut(G)$, then
$$
\varphi(b^s a^t c^r) = b^{v_\alpha t + v_\beta s} a^{u_\alpha t + u_\beta s} c^{ \Delta r + w_\alpha t + w_\beta s + \frac{1}{2} u_\alpha v_\alpha t(t-1) + \frac{1}{2} u_\beta v_\beta s(s-1) + u_\beta v_\alpha t s},
$$
where $0 \leq t, s, r \leq p-1$.

We define $\varphi_\gamma$ on $\mathbb{Z}_{p^2}$ for $G = E_3$ in this way: $\varphi_\gamma(x, y) = \varphi_\gamma(x \ \mathrm{rem}\ p, y \ \mathrm{rem}\ p)$ for $x, y \in \mathbb{Z}_{p^2}$.We note that such an extension modulo $p$ coincides with the right side of equality \eqref{phi_gamma} considered in the case $t, s\in \mathbb{Z}_{p^2}$.

Let us define a map $\varphi$ on~$G = E_3$ by formulas \eqref{define_automorphism}, \eqref{phi_alpha_beta} and \eqref{phi_gamma}, assuming that $t, s, r \in \mathbb{Z}_{p^2}$.
We show that
$\varphi(b^y a^x c^z) = \varphi(b^{y \ \mathrm{rem} \ p} a^{x \ \mathrm{rem} \ p} c^{(z + [x/p] + [y/p]) \ \mathrm{rem} \ p})$ for any $x, y, z\in \mathbb{Z}_{p^2}$,
using the fact that $u_\alpha + v_\alpha = u_\beta + v_\beta = \Delta \pmod{p}$:
\begin{multline*}
\varphi(b^y a^x c^z) = \varphi(b^y a^x) \varphi(c^z )= b^{v_\alpha x + v_\beta y} a^{u_\alpha x + u_\beta y} c^{\varphi_\gamma(x, y, 0) + kz} \\
= b^{\varphi_\beta(x \ \mathrm{rem}\ p, y \ \mathrm{rem}\ p)} b^{v_\alpha[x/p]p + v_\beta[y/p]p} a^{\varphi_\alpha(x \ \mathrm{rem}\ p, y \ \mathrm{rem}\ p)} a^{u_\alpha[x/p]p + u_\beta[y/p]p} c^{\varphi_\gamma(x \ \mathrm{rem} \ p, y \ \mathrm{rem} \ p, 0) + \Delta z} \\
= b^{\varphi_\beta(x \ \mathrm{rem}\ p, y \ \mathrm{rem}\ p)} a^{\varphi_\alpha(x \ \mathrm{rem}\ p, y \ \mathrm{rem}\ p)} c^{(v_\alpha + u_\alpha)[x/p] + (v_\beta + u_\beta)[y/p] + \varphi_\gamma(x \ \mathrm{rem} \ p, y \ \mathrm{rem} \ p, 0) + \Delta z} \\
= b^{\varphi_\beta(x \ \mathrm{rem}\ p, y \ \mathrm{rem}\ p)} a^{\varphi_\alpha(x \ \mathrm{rem}\ p, y \ \mathrm{rem}\ p)} c^{\varphi_\gamma(x \ \mathrm{rem} \ p, y \ \mathrm{rem} \ p, 0) + \Delta (z + [x/p] + [y/p])} \\
= \varphi(b^{y \ \mathrm{rem} \ p} a^{x \ \mathrm{rem} \ p} c^{(z + [x/p] + [y/p]) \ \mathrm{rem} \ p}).
\end{multline*}

For $G = H_3$, a similar remark is obvious since $a^p = b^p = e$.

Let $G = E_3$ and $\varphi \in \Aut(G)$ with parameters $u_\alpha, u_\beta, v_\alpha, v_\beta, w_\alpha, w_\beta \in \mathbb{Z}_{p^2}$ and $\Delta \in \mathbb{Z}_{p^2}^{*}$.In this case $u_\beta = u_\alpha - 1 \pmod{p}$, $\Delta = u_\alpha + v_\alpha = u_\beta + v_\beta \pmod{p}$. We obtain
$$
f_u = u_\beta - u_\alpha + 1 \in p \mathbb{Z}_{p^2}, \quad 
f_\alpha = \Delta - u_\alpha - v_\alpha - f_u \in p \mathbb{Z}_{p^2}, \quad 
f_\beta = \Delta - u_\beta - v_\beta \in p \mathbb{Z}_{p^2},
$$
we define
\begin{gather*}
u_\alpha' = u_\alpha + f_u, \quad
u_\beta' = u_\beta, \quad
v_\alpha' = v_\alpha + f_\alpha, \\
v_\beta' = v_\beta + f_\beta, \quad
w_\alpha' = w_\alpha - [(f_u + f_\alpha)/p], \quad
w_\beta' = w_\beta - [f_\beta/p].
\end{gather*}
The map defined by formulas \eqref{define_automorphism}, \eqref{phi_alpha_beta} and \eqref{phi_gamma} with parameters $u_\alpha'$, $u_\beta'$, $v_\alpha'$, $v_\beta'$, $w_\alpha'$, $w_\beta'$, $\Delta$ coincides with $\varphi$, and relations on the parameters hold in $\mathbb{Z}_{p^2}$:
$u_\beta' = u_\alpha' - 1$, $u_\alpha' + v_\alpha' = u_\beta + v_\beta' = \Delta$. 
Thus, we can assume that the parameters of the automorphisms~$\varphi$ of the group~$G$ satisfy these three relations on $\mathbb{Z}_{p^2}$.

Now we consider the inverse to this automorphism $\varphi$:
\begin{gather*}
\varphi^{-1}(b^s a^t c^r) = b^{v_\alpha' t + v_\beta' s} a^{u_\alpha' t + u_\beta' s} c^{ \Delta' r + w_\alpha' t + w_\beta' s + \frac{1}{2} u_\alpha' v_\alpha' t(t-1) + \frac{1}{2} u_\beta' v_\beta' s(s-1) + u_\beta' v_\alpha' t s},
\end{gather*}
where
\begin{equation} \label{op_automorphism}
\begin{gathered}
u_\alpha' = \frac{v_\beta}{\Delta}, \quad
u_\beta' = \frac{-u_\beta}{\Delta}, \quad
v_\alpha' = \frac{-v_\alpha}{\Delta}, \quad
v_\beta' = \frac{u_\alpha}{\Delta}, \quad
\Delta' = \frac{1}{\Delta}, \\ 
w_\alpha' = \frac{1}{\Delta^2}\left(-w_\alpha v_\beta + w_\beta v_\alpha + \frac{u_\alpha v_\alpha v_\beta}{2} - \frac{u_\beta v_\alpha v_\beta}{2} - \frac{v_\alpha v_\beta}{2}\right), \\ 
w_\beta' = \frac{1}{\Delta^2}\left(w_\alpha u_\beta - w_\beta u_\alpha -\frac{u_\alpha u_\beta v_\alpha}{2} + \frac{u_\alpha u_\beta v_\beta}{2} - \frac{u_\alpha u_\beta}{2}\right),
\end{gathered}
\end{equation}
where $u_\alpha, u_\beta, v_\alpha, v_\beta, w_\alpha, w_\beta \in \mathbb{Z}_p$, $\Delta = u_\alpha v_\beta - u_\beta v_\alpha \in \mathbb{Z}_{p}^{*}$ for $G = H_3$ and $u_\alpha, u_\beta, v_\alpha, v_\beta, w_\alpha, \\
w_\beta \in \mathbb{Z}_{p^2}$, $\Delta \in \mathbb{Z}_{p^2}^{*}$ and $u_\beta = u_\alpha - 1$, $\Delta = u_\alpha + v_\alpha = u_\beta + v_\beta$ for $G = E_3$.

The fact that this automorphism indeed is the inverse is shown in Appendix A.

Conjugating a Rota---Baxter operator by automorphisms, we also obtain a Rota---Baxter operator. We denote $B^{(\varphi)} = \varphi^{-1} \circ B \circ \varphi$.
Let $\varphi \in \Aut(G)$, we will use the following notations: 
$(l_\alpha^{(\varphi)}, l_\beta^{(\varphi)}, m_\alpha^{(\varphi)}, m_\beta^{(\varphi)}, x_\alpha^{(\varphi)}, x_\beta^{(\varphi)}, k^{(\varphi)}) = \chi (B^{(\varphi)})$ for $G = H_3$ and $[l_\alpha^{(\varphi)}, l_\beta^{(\varphi)}, m_\alpha^{(\varphi)}, m_\beta^{(\varphi)}, x_\alpha^{(\varphi)}, x_\beta^{(\varphi)}, k^{(\varphi)}] = \chi_E (B^{(\varphi)})$ for $G = E_3$.

\section{Invariance of the center of an extraspecial group}
Let $G$ be an extraspecial group of order $p^{1+2n}$.  Then a following presentation exists:
\begin{multline*}
G = \langle a_1, \dots, a_n, b_1, \dots, b_n, c \mid [a_i, b_i] = c, \, a_m^p = b_m^p = c^p = [a_i, a_j] = [b_i, b_j] \\
= [a_k, b_l] = [a_i, c] = [b_i, c] = e \rangle, 
\end{multline*}
where $k \neq l$ and either $1 \leq m \leq n$, or $1 \leq m \leq n-1$, $a_n^{p} = b_n^{p} = c$. 
We will denote $Z(G) = \langle c \rangle$ as~$Z$.

Further, we will use arguments from the paper \cite{Gorenstein} to prove item b) of Lemma 1 and to prove Lemma 2.

Since $G/Z \cong \bigoplus_{i =1}^{2n} \mathbb{Z}_p$, we can consider $G/Z$ as vector space over $\mathbb{Z}_p$.
We denote $V = G/Z$. We define a symplectic form $f\colon V \otimes V \to \mathbb{Z}_p$ as $f(u, v) = [u, v] \in Z \cong \mathbb{Z}_p$. 

\textbf{Lemma 1}.
Let $H$ be subgroup~$G$. Then

a) if $H$ is nonabelian, then $Z\subseteq H$;

b) if $|H| \geq p^{n+2}$, then $H$ is nonabelian;

c) if $H$ is nonabelian, $K \trianglelefteq H$ and $c \notin K$, then $K \subseteq Z(H)$;

d) if $|H| \geq p^{n+1}$, then $Z\subseteq H$.

{\sc Proof}.
a) If $h_1, h_2 \in H$ do not commute, then $\langle [h_1, h_2] \rangle = Z$.

b) Let $A$ be a maximal abelian subgroup. Then $Z \subseteq A$, holds, otherwise $\langle A, c \rangle$ is abelian and~$A\subsetneq \langle A, c \rangle$.
Thus, $U = A/Z$ is the subspace in $V$. We have $f(U, U) = 0$, since $A$ is abelian.
It is known that the dimension of an isotropic subspace in a $2n$-dimensional symplectic space does not exceed $n$.
Consequently, $\dim(U) \leq n$ holds and, thus, $|A| \leq p^{n+1}$ holds, the statement of item b) is proved..

c) Assume that $K\not\subseteq Z(H)$, then we can find $g \in H \setminus Z(H)$ and $h \in K$, such that $gh\neq hg$. 
Consequently, $1 \leq k \leq p-1$ exists such that $c^k = g^{-1} h^{-1} g\cdot h\in K$. We obtained a contradiction.

d) If $H$ is nonabelian, then $Z\subseteq H$ holds by item a). If $H$ is abelian, then $\langle H, c \rangle$ also is abelian and $|\langle H, c \rangle| \geq p^{n+2}$ holds, we obtained a contradiction with b).
\hfill $\square$

\textbf{Lemma 2}.
Let $H \leq G$ and $|H| = p^{n+1+k}$, where $1 \leq k \leq n$. Then $|Z(H)| \leq p^{n+1-k}$ and $c \in Z(H)$ hold.

{\sc Proof}.
An inclusion $Z \subseteq Z(H)$ exists by Lemma 1a), b), hence, $U = H/Z$ is subspace in~$V$ and $\dim\,U = n+k$.

We denote by $U^\perp$ the orthogonal complement to $U$ in $V$ with respect to $f$. 
Then $Z(H)/Z = U \cap U^{\perp}$ holds since $f(Z(H)/Z, H/Z) = 0$ and $Z(H) \subseteq H$ hold. 
We have $\dim(U^\perp) = \dim(V) - \dim(U) = n-k$. 
Therefore $\dim(U \cap U^{\perp}) \leq n-k$ holds and we obtain $|Z(H)| \leq |Z| p^{n-k} = p^{n+1-k}$.~\hfill $\square$

\textbf{Lemma 3}.
Let br $G = H_3$ or $G = E_3$ and $B$ be an RB-operator on $G$. Then $B(Z) \subseteq Z$.

{\sc Proof}.
Let us denote $H = \text{Im}(B)$, $\widetilde{H} = \text{Im}(\widetilde{B})$.
Up to the action\ \ $\widetilde{}$\, we can assume that $|H| \leq |\widetilde{H}|$.

{\sc Case 1}: $|H| = |\widetilde{H}| = p^3$. 
Then $B(Z) \subseteq Z$ \cite[Lemma 2]{GG}.

{\sc Case 2}:
Let $|H| \leq p^2$, $|\widetilde{H}| = p^3$. Consequently, $\ker(B) \trianglelefteq G$, $|\ker(B)| \geq p$, and therefore $c \in \ker(B)$ by Lemma 1v).

{\sc Case 3}:
Let $|H| = p^2$, $|\widetilde{H}| = p^2$. In this case, $Z \subseteq H, \widetilde{H}$ by Lemma 1d) and $|H \cap \widetilde{H}| = p$ by Proposition 4b). Hence $H \cap \widetilde{H} = Z$. Then $B(Z) \subseteq B(\widetilde{H}) = H \cap \widetilde{H} = Z$.

{\sc Case 4}:
Let $|H| = p$, $|\widetilde{H}| = p^2$. Then $Z \subseteq \widetilde{H}$ by Lemma 1g) and $H \cap \widetilde{H} = \{ e \}$ by Proposition 4b). Consequently, $B(Z) \subseteq B(\widetilde{H}) = H \cap \widetilde{H} = \{ e \}$.
\hfill $\square$

\textbf{Theorem 1}.
Let $B$ be a Rota---Baxter operator on~$G$. Then $B(Z) \subseteq Z$.

{\sc Proof}.
Let us prove the statement by induction on~$n$.
For $n=1$, this statement is Lemma~3. Assume that $n \geq 2$ and the statement is proved for all $k<n$.
Below we will apply the induction hypothesis only for Case 5.

We denote $H = \Imm(B)$, $\widetilde{H} = \Imm(\widetilde{B})$
and we may assume that
$|H| \leq |\widetilde{H}|$.

{\sc Case 1}: $|H| = |\widetilde{H}| = p^{2n+1}$. 
Then $B(Z) \subseteq Z$ \cite[lemma 2]{GG}.

{\sc Case 2}: $|H| \leq p^{n}$. 
Then $|\ker(B)| \geq p^{n + 1}$ and $Z \subseteq \ker(B)$ by lemma 1d).

{\sc Case 3}: $|H| = |\widetilde{H}| = p^{n+1}$. 
Then $Z \subseteq H, \widetilde{H}$ by lemma 1 d), moreover, $|H \cap \widetilde{H}| = p$ due to $H \cdot \widetilde{H} = G$. Then $H \cap \widetilde{H} = Z$ and we obtain $\Imm(B \widetilde{B}) = Z$. We have $B(\widetilde{H}) = Z$, hence $B(Z) \subseteq Z$, since $Z \subseteq \widetilde{H}$.

{\sc Case 4}: $|H| < p^{n+1+k}$, $|\widetilde{H}| = p^{n+1+k}$, where $1 \leq k \leq n$. 
Then $|\ker(B)| \geq p^{n+1-k}$, an we have $|Z(\widetilde{H})| \leq p^{n+1-k}$ by lemma 2, and $c\in Z(\widetilde{H})$. 
By item d) of Lemma 1, the subgroup~$\widetilde{H}$ is non-abelian.
Assume that $c\not \in \ker (B)$.
Then, applying Lemma 1v), we conclude that $Z(\widetilde{H}) = \ker(B)$, therefore $c \in \ker(B)$.

{\sc Case 5}: $|H| = |\widetilde{H}| = p^{n+1+k}$, where $1 \leq k \leq n-1$. 
Then $|\ker(B)| = |\ker(\widetilde{B})| = p^{n - k}$ and $|Z(H)|, |Z(\widetilde{H})| \leq p^{n+1-k}$ by lemma 2. 
We have $Z(H) \cap Z(\widetilde{H}) =Z$ by lemma~2 and by $H \cdot \widetilde{H} = G$. 
If $c \in \ker(B) \cup \ker(\widetilde{B})$, then $B(Z) \subseteq Z$. Let $c \notin \ker(B), \ker(\widetilde{B})$. Then we obtain
$$
\ker(B)\subsetneq Z(\widetilde{H}), \quad
\ker(\widetilde{B})\subsetneq Z(H), \quad
|Z(H)| = |Z(\widetilde{H})| = p^{n+1-k}
$$ 
by Lemma 1c), and since $c \in Z(H)\cap Z(\widetilde{H})$.

The operator $B|_{\widetilde{H}}$ is well-defined. Let us consider $K = \widetilde{H} / \ker(B)$. Then the Rota---Baxter operator $B_K \colon K \to K$ is naturally defined on the quotient group. 
We specify some properties of~$K$. 
The group $K$ is non-abelian and $Z(K) = \langle c\cdot\ker B \rangle$. 
For any $g, h \in \widetilde{H}$, we have $[g\cdot \ker B, h\cdot \ker B] \in Z(K)$ and 
$(g\cdot \ker B)^p \in Z(K)$. 
Then $K /Z(K) \cong \mathbb{Z}_p^{2k}$. 
We conclude that $K$ is an extraspecial group of order $2k + 1 < 2n + 1$. 
Applying the induction hypothesis, we obtain that $B_K(c\cdot\ker B) = (c\cdot\ker B)^r$ for some $0\leq r\leq p-1$. 
Consequently, $B(c) \in Z(\widetilde{H})$. 

Considering the Rota---Baxter operator $\widetilde{B}$ on the group $L = H / \ker(\widetilde{B})$ and repeating the arguments above, we deduce that $\widetilde{B}(c) \in Z(H)$ and, thus, $\widetilde{B}(c^{-1}) \in Z(H)$ since $\widetilde{B}(c)\widetilde{B}(c^{-1}) = e$.
On the one hand, $B(c) \in Z(\widetilde{H})$. On the other hand, 
$B(c)  = c^{-1} \widetilde{B}(c^{-1}) \in Z(H)$. 
Hence $B(c) \in Z(H) \cap Z(\widetilde{H}) = Z$.
\hfill $\square$

\section{Rota---Baxter Operators on $H_3$ and $E_3$}
\subsection{General Form of Rota---Baxter Operators on $H_3$ and $E_3$}

In this and the next section, we assume that $G = H_3$ or $G = E_3$.

Consider the subset $\{0,\ldots,p-1\}$ of the set $\mathbb{Z}_{p^2}$.
Any element $g \in G$ is uniquely representable in the form $g = b^s a^t c^r$, where $0 \leq t,s,r \leq p-1$.

Fix an RB-operator $B$ on $G$. Then we express the action of $B$ on the elements as follows:
\begin{equation} \label{B-Main}
B(b^s a^t c^r) = b^{B_\beta(t, s)} a^{B_\alpha(t, s)} c^{B_\gamma(t, s, r)}.
\end{equation}
If $G = H_3$, we assume that 
$B_\alpha, B_\beta \colon 
\mathbb{Z}_{p} \times \mathbb{Z}_{p} \to \mathbb{Z}_{p}$, $B_\gamma \colon \mathbb{Z}_{p} \times \mathbb{Z}_{p} \times \mathbb{Z}_{p} \to \mathbb{Z}_{p}$,
$k\in \mathbb{Z}_{p}$. If $G = E_3$, then we take $B_\alpha, B_\beta \colon 
\{0,\ldots,p-1\} \times \{0,\ldots,p-1\} \to \mathbb{Z}_{p^2}$, $B_\gamma \colon \{0,\ldots,p-1\} \times \{0,\ldots,p-1\} \times \{0,\ldots,p-1\} \to \mathbb{Z}_{p}$,
$k\in \mathbb{Z}_{p^2}$.

Since
$$
B(g) = B(b^s a^t c^r) = B(c^r B(c^r) b^s a^t B(c^r)^{-1}) = B(c^r) B(b^s a^t),
$$
the equality $B_\gamma(t, s, r) = B_\gamma(t, s, 0) + B_\gamma(0, 0, r)$ holds.

Due to the fact that $B(Z) \subseteq Z$, we obtain that
\begin{equation}\label{B_alpha_beta}
B_\alpha(t, s) = l_\alpha t + l_\beta s, \quad 
B_\beta(t, s) = m_\alpha t + m_\beta s, \quad
B_\gamma(0, 0, r) = kr,
\end{equation}
where $l_\alpha, l_\beta, m_\alpha, m_\beta, k \in \mathbb{Z}_p$ if $G = H_3$, and $l_\alpha, l_\beta, m_\alpha, m_\beta, k \in \mathbb{Z}_{p^2}$ if $G = E_3$.
We extend the mapping $B_\alpha, B_\beta\colon \mathbb{Z}_{p^2} \times \mathbb{Z}_{p^2} \to \mathbb{Z}_{p^2}$
by formula~\eqref{B_alpha_beta}, where $t,s\in \mathbb{Z}_{p^2}$ in the case when $G = E_3$.

Consider $0 \leq t_1, t_2, s_1, s_2 \leq p-1$; then in the case $G = E_3$, we obtain the following relation:
\begin{gather*}
B_\alpha(t_1, s_1) + B_\alpha(t_2, s_2) = l_\alpha (t_1 + t_2) + l_\beta (s_1 + s_2) = l_\alpha(t_1 + t_2 \ \mathrm{rem} \ p) \\
+ l_\beta (s_1 + s_2 \ \mathrm{rem} \ p) + l_\alpha p [(t_1 + t_2) / p] + l_\beta p [(s_1 + s_2) / p] \allowdisplaybreaks \\ 
= B_\alpha(t_1 + t_2 \ \mathrm{rem} \ p, s_1 + s_2 \ \mathrm{rem} \ p) + l_\alpha p [(t_1 + t_2) / p] + l_\beta p [(s_1 + s_2) / p].
\end{gather*}
A similar formula holds for \(B_{\beta }\). We determine the form of \(B_{\gamma }\) and the necessary parameter restrictions by applying the Rota–Baxter operator definition to the elements \(b^{s_1} a^{t_1}\) and \(b^{s_2} a^{t_2}\).
\begin{multline*}
b^{B_\beta((t_1 + t_2)\ \mathrm{rem}\ p, (s_1 + s_2)\ \mathrm{rem}\ p)} a^{B_\alpha((t_1 + t_2)\ \mathrm{rem}\ p, (s_1 + s_2)\ \mathrm{rem}\ p)} \\
\times c^{ B_\alpha(t_1, s_1) B_\beta(t_2, s_2) + B_\gamma(t_1, s_1) + B_\gamma(t_2, s_2) + E(t_1 + t_2, s_1 + s_2)} \\
= b^{B_\beta(t_1, s_1) + B_\beta(t_2, s_2)} a^{B_\alpha(t_1, s_1) + B_\alpha(t_2, s_2)} c^{ B_\alpha(t_1, s_1) B_\beta(t_2, s_2) + B_\gamma(t_1, s_1, 0) + B_\gamma(t_2, s_2, 0)} \\
 = B(b^{s_1} a^{t_1}) B(b^{s_2} a^{t_2}) = B( b^{s_1} a^{t_1} B(b^{s_1} a^{t_1}) b^{s_2} a^{t_2} B(b^{s_1} a^{t_1})^{-1} ) \\
 = B(b^{s_1} a^{t_1} b^{B_\beta(t_1, s_1)} a^{B_\alpha(t_1, s_1)} b^{s_2} a^{t_2} a^{-B_\alpha(t_1, s_1)} b^{-B_\beta(t_1, s_1)}) \\
 = B(b^{(s_1 + s_2)\ \mathrm{rem}\ p} a^{(t_1 + t_2)\ \mathrm{rem}\ p} c^{t_1 s_2 + s_2 B_\alpha(t_1, s_1) - t_2 B_\beta(t_1, s_1) + R(t_1 + t_2, s_1 + s_2)}) \\ 
 = b^{B_\beta((t_1 + t_2)\ \mathrm{rem}\ p, (s_1 + s_2)\ \mathrm{rem}\ p)} a^{B_\alpha((t_1 + t_2)\ \mathrm{rem}\ p, (s_1 + s_2)\ \mathrm{rem}\ p)} \\
 \times c^{B_\gamma(t_1 + t_2 \ \mathrm{rem}\ p, s_1 + s_2 \ \mathrm{rem}\ p, 0) + k( t_1 s_2 + s_2 B_\alpha(t_1, s_1) - t_2 B_\beta(t_1, s_1) + R(t_1 + t_2, s_1 + s_2) )},
\end{multline*}
where $R, E \colon \mathbb{Z}_{p^2} \times \mathbb{Z}_{p^2} \to \mathbb{Z}_{p^2}$, 
$$
R(t, s) = [{t}/{p}] + [{s}/{p}], \quad
E(t, s) =  [{t}/{p}] (l_\alpha + m_\alpha) + [{s}/{p}] (l_\beta + m_\beta)
$$ 
for $G = E_3$ and $R(t, s) = E(t, s) = 0$ for $G = H_3$.

Then we obtain the relation modulo $p$:
\begin{multline*}
B_\gamma(t_1, s_1, 0) + B_\gamma(t_2, s_2, 0) + (l_\alpha m_\alpha t_1 t_2 + l_\alpha m_\beta t_1 s_2 + l_\beta m_\alpha t_2 s_1 + l_\beta m_\beta s_1 s_2) + E(t_1 + t_2, s_1 + s_2) \\
= B_\gamma(t_1, s_1, 0) + B_\gamma(t_2, s_2, 0) + B_\alpha(t_1, s_1) B_\beta(t_2, s_2) + E(t_1 + t_2, s_1 + s_2) \\
= B_\gamma(t_1 + t_2 \ \mathrm{rem}\ p, s_1 + s_2 \ \mathrm{rem}\ p, 0) + k( t_1 s_2 + s_2 B_\alpha(t_1, s_1) - t_2 B_\beta(t_1, s_1) + R(t_1 + t_2, s_1 + s_2) ) \\
= B_\gamma(t_1 + t_2 \ \mathrm{rem}\ p, s_1 + s_2 \ \mathrm{rem}\ p, 0) + k( t_1 s_2 + l_\alpha t_1 s_2 + l_\beta s_1 s_2 - m_\alpha t_1 t_2 - m_\beta t_2 s_1 \\
 + R(t_1 + t_2, s_1 + s_2) ).
\end{multline*}

In particular, the equalities below modulo $p$ follow from this relation:
\begin{gather*}
B_\gamma(t, 0, 0) + B_\gamma(0, s, 0) + l_\alpha m_\beta t s = B_\gamma(t, s, 0) + k(l_\alpha + 1) t s + k R(t, s) - E(t, s), \allowdisplaybreaks \\
B_\gamma(0, s, 0) + B_\gamma(t, 0, 0) + l_\beta m_\alpha t s = B_\gamma(t, s, 0) - k m_\beta t s + k R(t, s) - E(t, s),
\end{gather*}

\vspace{-1.4cm}

\begin{multline*}
B_\gamma(t_1, 0, 0) + B_\gamma(t_2, 0, 0) + l_\alpha m_\alpha t_1 t_2 \\
= B_\gamma(t_1 + t_2 \ \mathrm{rem}\ p, 0, 0) - k m_\alpha t_1 t_2 + k R(t_1 + t_2, 0) - E(t_1 + t_2, 0), 
\end{multline*}

\vspace{-1.4cm}

\begin{multline*}
B_\gamma(0, s_1, 0) + B_\gamma(0, s_2, 0) + l_\beta m_\beta s_1 s_2 \\
 = B_\gamma(0, s_1 + s_2 \ \mathrm{rem}\ p, 0) + k l_\beta s_1 s_2 + k R(0, s_1 + s_2) - E(0,  s_1 + s_2), 
\end{multline*}

\vspace{-0.9cm}

$$
l_\alpha m_\beta - k l_\alpha = l_\beta m_\alpha + k m_\beta + k.
$$

They are equivalent to the following four formulas modulo $p$:
\begin{equation}\label{pre_B_gamma_a_b} 
B_\gamma(t, s,0) = B_\gamma(t, 0,0) + B_\gamma(0, s,0) + (l_\beta m_\alpha + k m_\beta) t s - k R(t, s) + E(t, s), 
\end{equation}

\vspace{-1.4cm}

\begin{multline}
B_\gamma(t_1 + t_2 \ \mathrm{rem}\ p, 0,0) \\
 = B_\gamma(t_1, 0,0) + B_\gamma(t_2, 0,0) + m_\alpha(l_\alpha + k) t_1 t_2 - k R(t_1 + t_2, 0) + E(t_1 + t_2, 0), \label{pre_B_gamma_a} 
\end{multline} 

\vspace{-1.4cm}

\begin{multline}
B_\gamma(0, s_1 + s_2 \ \mathrm{rem}\ p,0) \\
= B_\gamma(0, s_1,0) + B_\gamma(0, s_2,0) + l_\beta(m_\beta - k) s_1 s_2 - k R(0, s_1 + s_2) + E(0,  s_1 + s_2), \label{pre_B_gamma_b}
\end{multline} 

\vspace{-0.9cm}

\begin{equation}\label{RB_condition}
l_\alpha m_\beta - k l_\alpha = l_\beta m_\alpha + k m_\beta + k. 
\end{equation}

We note that $R(t, s) = E(t, s) = 0$ for $0 \leq t, s \leq p-1$.

\textbf{Lemma 4}.
Let $0 \leq t, s \leq p-1$. Then the formulas hold modulo $p$
\begin{gather}
B_\gamma(t, 0, 0) = B_\gamma(1, 0, 0) t + \frac{t(t-1)}{2} m_\alpha(l_\alpha + k), \label{B_gamma_a} \\
B_\gamma(0, s, 0) = B_\gamma(0, 1, 0) s + \frac{s(s-1)}{2} l_\beta(m_\beta - k) \label{B_gamma_b}.
\end{gather}

{\sc Proof}.
We prove formula~\eqref{B_gamma_a}.
We carry out the proof by induction on~$t$ from $0$ to $p-1$.
The induction base for $t = 0$ holds since $B_\gamma(0, 0, 0) = 0$.
Assume that \eqref{B_gamma_a} is true for $t < p-1$; we prove it for $t+1$. 
We use formula~\eqref{pre_B_gamma_a} and obtain:
\begin{multline*}
B_\gamma(t+1, 0, 0) = B_\gamma(t, 0, 0) + B_\gamma(1, 0, 0) + m_\alpha (l_\alpha + k)t\allowdisplaybreaks \\ 
= B_\gamma(1, 0, 0) t + \frac{t(t-1)}{2} m_\alpha (l_\alpha + k) + B_\gamma(1, 0, 0) + l_\alpha m_\alpha t + k m_\alpha t\\ 
= B_\gamma(1, 0, 0) (t+1) + \frac{(t+1)t}{2} m_\alpha(l_\alpha + k).
\end{multline*}
Formula \eqref{B_gamma_a} is proved.
The proof of ~\eqref{B_gamma_b} is carried out similarly.
\hfill $\square$

We denote $x_\alpha = B_\gamma(1, 0, 0)$, $x_\beta = B_\gamma(0, 1, 0)$.
We obtain
\begin{multline} \label{B_gamma}
B_\gamma(t, s, r) = x_\alpha t + x_\beta s + (l_\beta m_\alpha + k m_\beta )ts + m_\alpha(l_\alpha + k)\frac{t(t-1)}{2} \\
+ l_\beta(m_\beta - k)\frac{s(s-1)}{2} + kr,
\end{multline}
where $0 \leq t,s,r <p$, by \eqref{pre_B_gamma_a_b}, \eqref{B_gamma_a} and \eqref{B_gamma_b}.

Thus, $B_\gamma$ is defined on $\{ 0, \ldots, p-1 \}$ in the case $G = E_3$; 
we determine $B_\gamma$ on the whole $\mathbb{Z}_{p^2}$ by the formula
$B_\gamma(x, y) = B_\gamma(x \ \mathrm{rem} \ p, y \ \mathrm{rem} \ p)$ for $x, y \in \mathbb{Z}_{p^2}$.
We note that such an extension modulo $p$ coincides with the right-hand side of equality \eqref{B_gamma} considered in the case when $t, s\in \mathbb{Z}_{p^2}$.

For $G = H_3$, the mappings $B_\alpha, B_\beta, B_\gamma$ define a Rota---Baxter operator if \eqref{RB_condition} holds. 

Consider the case $G = E_3$. In this case, we additionally require the fulfillment of the equality $E(t, s) = kR(t, s) \pmod{p}$ as a necessary and sufficient condition because $B_\gamma(t_1 + t_2 \ \mathrm{rem}\ p, s_1 + s_2 \ \mathrm{rem}\ p, 0) = B_\gamma(t_1 + t_2, s_1 + s_2)$ modulo $p$; consequently, $k = l_\alpha + m_\alpha = l_\beta + m_\beta \pmod{p}$. Substituting these identities into \eqref{RB_condition}, we obtain the condition $k^2 + k = 0 \pmod{p}$, that is, $k = 0, -1 \pmod{p}$. 
Summarizing the above, we obtain
\begin{equation} \label{RB_condition_E_3}
k = l_\alpha + m_\alpha = l_\beta + m_\beta \pmod{p}; \quad k = 0, -1 \pmod{p}.
\end{equation}

Thus, we obtained

\textbf{Proposition 5}.
Let $B$ be Rota---Baxter operator on $G = H_3, E_3$. Then maps $B_\alpha$, $B_\beta$, $B_\gamma$ have following forms:
\begin{gather*}
B_\alpha(t, s) = l_\alpha t + l_\beta s, \quad 
B_\beta(t, s) = m_\alpha t + m_\beta s, \\
B_\gamma(t, s, r) = x_\alpha t + x_\beta s + (l_\beta m_\alpha + k m_\beta )ts + m_\alpha(l_\alpha + k)\frac{t(t-1)}{2} + l_\beta(m_\beta - k)\frac{s(s-1)}{2} + kr,
\end{gather*}
where the following equalities are satisfied:

$l_\alpha m_\beta - k l_\alpha = l_\beta m_\alpha + k m_\beta + k  \pmod{p}$ for $G = H_3$;

$k = l_\alpha + m_\alpha = l_\beta + m_\beta \pmod{p}, \quad k = 0, -1 \pmod{p}$ for $G = E_3$.

\textbf{Remark 2}.
In the case $G = E_3$, choosing different parameters $l_\alpha, l_\beta, m_\alpha, m_\beta, x_\alpha, x_\beta, k$ from $\mathbb{Z}_{p^2}$,
we can obtain the same Rota---Baxter operator.

Let us consider the mapping $B$ on~$G = E_3$ defined by formulas~\eqref{B-Main},~\eqref{B_alpha_beta}, \eqref{B_gamma}, assuming that $t, s, r \in \mathbb{Z}_{p^2}$.
Now we show that for any $x,y,z\in \mathbb{Z}_{p^2}$, the equality
$B(b^y a^x c^z) = B(b^{y \ \mathrm{rem} \ p} a^{x \ \mathrm{rem} \ p} c^{(z + [x/p] + [y/p]) \ \mathrm{rem} \ p})$ holds, using the identities $k = l_\alpha + m_\alpha = l_\beta + m_\beta \pmod{p}$:
\begin{multline*}
B(b^y a^x c^z) = B(b^y a^x) B(c^z)= B(b^y a^x) c^{k z} = b^{m_\alpha x + m_\beta y} a^{l_\alpha x + l_\beta y} c^{B_\gamma(x, y, 0) + kz} \\
= b^{B_\beta(x \ \mathrm{rem}\, p, y \ \mathrm{rem}\, p)} b^{m_\alpha[x/p]p + m_\beta[y/p]p} a^{B_\alpha(x \ \mathrm{rem}\, p, y \ \mathrm{rem}\, p)} a^{l_\alpha[x/p]p + l_\beta[y/p]p} c^{B_\gamma(x \ \mathrm{rem} \, p, y \ \mathrm{rem} \, p, 0) + kz} \\
= b^{B_\beta(x \ \mathrm{rem}\ p, y \ \mathrm{rem}\ p)} a^{B_\alpha(x \ \mathrm{rem}\ p, y \ \mathrm{rem}\ p)} c^{(m_\alpha + l_\alpha)[x/p] + (m_\beta + l_\beta)[y/p] + B_\gamma(x \ \mathrm{rem} \ p, y \ \mathrm{rem} \ p, 0) + kz} \\
= b^{B_\beta(x \ \mathrm{rem}\ p, y \ \mathrm{rem}\ p)} a^{B_\alpha(x \ \mathrm{rem}\ p, y \ \mathrm{rem}\ p)} c^{B_\gamma(x \ \mathrm{rem} \ p, y \ \mathrm{rem} \ p, 0) + k(z + [x/p] + [y/p])} \\
= B(b^{y \ \mathrm{rem} \ p} a^{x \ \mathrm{rem} \ p} c^{(z + [x/p] + [y/p]) \ \mathrm{rem} \ p}).
\end{multline*}

For $G = H_3$, a similar remark is obvious since $a^p = b^p = e$.

Let $\mathrm{RB}(G)$ be the set of Rota---Baxter operators on $G$.
For $G = H_3$, we define a mapping $\chi \colon \mathrm{RB}(G) \to \mathbb{Z}_p^7$ as follows: $\chi(B) = (l_\alpha, l_\beta, m_\alpha,m_\beta$, $x_\alpha, x_\beta, k)$,
where $B$ is a Rota---Baxter operator defined by formulas \eqref{B_alpha_beta} and \eqref{B_gamma}.
We will denote the $i$-th coordinate of $\chi(B)$ as $\chi(B)_i$, where $1 \leq i \leq 7$.

Now we define an analogue of the mapping $\chi$ for $G = E_3$. 
Let $[A]$ denote the collection of equivalence classes on the set of Rota---Baxter parameters $A \subseteq \mathbb{Z}_{p^2}^7$. Here, the tuples $(l_\alpha, l_\beta, m_\alpha, m_\beta$, $x_\alpha,x_\beta, k)$ and $(l_\alpha', l_\beta', m_\alpha', m_\beta', x_\alpha', x_\beta', k')$ are equivalent if these sets define the same Rota---Baxter operator by formulas \eqref{B_alpha_beta} and \eqref{B_gamma}. We define a mapping $\chi_E \colon RB(G) \to A$ as $\chi_E(B) = [l_\alpha, l_\beta, m_\alpha, m_\beta, x_\alpha, x_\beta]$ for an RB-operator $B$ with the corresponding parameters.

Let $G = E_3$.
We choose $B$ such that $\chi_E(B) = [l_\alpha, l_\beta, m_\alpha, m_\beta, x_\alpha, x_\beta, 0]$. 
We have $l_\alpha + m_\alpha = l_\beta + m_\beta = 0 \pmod{p}$ by \eqref{RB_condition_E_3}. Consider $f_\alpha = l_\alpha + m_\alpha$, $f_\beta = l_\beta + m_\beta$ and define
$$
l_\alpha' = l_\alpha - f_\alpha, \quad 
l_\beta' = l_\beta - f_\beta, \quad 
m_\alpha' = m_\alpha, \quad 
m_\beta' = m_\beta, \quad
x_\alpha' = x_\alpha + [f_\alpha/p], \quad 
x_\beta' = x_\beta + [f_\beta/p].
$$
Then $\chi_E(B) = [l_\alpha', l_\beta', m_\alpha', m_\beta', x_\alpha', x_\beta', 0]$ and $l_\alpha' + m_\alpha' = l_\beta' + m_\beta' = 0$ in $\mathbb{Z}_{p^2}$. Thus, we can assume that $l_\alpha + m_\alpha = l_\beta + m_\beta = 0$ in $\mathbb{Z}_{p^2}$.

\textbf{Lemma 5}.
a) Let $B$ be a Rota---Baxter operator on $H_3$ and let $\chi(B) = (l_\alpha, l_\beta, m_\alpha, m_\beta, $ $x_\alpha, x_\beta, k)$. Then $\chi(\widetilde{B}) = (-(l_\alpha+1), -l_\beta, -m_\alpha, -(m_\beta + 1), m_\alpha(l_\alpha + k + 1) - x_\alpha, l_\beta(m_\beta - k) - x_\beta, -(k+1) )$.

b) Let $B$ be a Rota---Baxter operator on $E_3$ and let $\chi_E(B) = [l_\alpha, l_\beta, m_\alpha, m_\beta,x_\alpha, x_\beta, k]$. Then $\chi_E(\widetilde{B}) = [-(l_\alpha+1), -l_\beta, -m_\alpha, -(m_\beta + 1), m_\alpha(l_\alpha + k + 1) - x_\alpha, l_\beta(m_\beta - k) - x_\beta,$ \\
$ -(k+1) ]$.

{\sc Proof}. 
We expand by definition:
\begin{multline*}
\widetilde{B}(b^s a^t c^r) = b^{-s} a^{-t} c^{ts-r} B(b^{-s} a^{-t} c^{ts-r}) = b^{-s} a^{-t} c^{ts-r} b^{-m_\alpha t - m_\beta s} a^{-l_\alpha t - l_\beta s} c^{B_\gamma(-t, -s, ts-r)} 
\\
= b^{-m_\alpha t - m_\beta (s + 1)} a^{-l_\alpha (t + 1) - l_\beta s} c^{t(m_\alpha t + m_\beta s) + ts - r + B_\gamma(-t, -s, ts-r)}.
\end{multline*}
Therefore, we obtain the following set of equalities modulo $p$:
\begin{multline*}
\widetilde{B}_\gamma(t, s, r) = t (m_\alpha t  + m_\beta s) + ts - r + B_\gamma(-t, -s, ts-r) = t (m_\alpha t  + m_\beta s) + ts - r \\
+ (l_\beta m_\alpha + k m_\beta ) ts + m_\alpha(l_\alpha + k) \frac{t (t + 1)}{2} + l_\beta(m_\beta - k) \frac{s (s + 1)}{2} - x_\alpha t - x_\beta s \\
+ k (ts - r) = ts - r + m_\alpha t^2 + m_\beta ts + ( (-l_\beta)(-m_\alpha) + (-k - 1)(-m_\beta - 1) ) ts \\
- (1 + m_\beta + k ) ts + m_\alpha(l_\alpha + k) \frac{t (t + 1)}{2} + l_\beta(m_\beta - k) \frac{s (s + 1)}{2} - x_\alpha t - x_\beta s + k (ts - r) = m_\alpha t^2 \\
+ ( (-l_\beta)(-m_\alpha) + (-k - 1)(-m_\beta - 1) ) ts + m_\alpha(l_\alpha + k) \frac{t (t + 1)}{2} + l_\beta(m_\beta - k) \frac{s (s + 1)}{2} - x_\alpha t \\
- x_\beta s - (k + 1)r = m_\alpha t^2 + ( (-l_\beta)(-m_\alpha) + (-k - 1)(-m_\beta - 1) ) ts - m_\alpha(-l_\alpha - 1 - k - 1) \frac{t (t - 1)}{2} \\
- m_\alpha t (t - 1) + m_\alpha(l_\alpha + k) t - l_\beta(-m_\beta - 1 - (-k - 1) ) \frac{s (s - 1)}{2} + l_\beta(m_\beta - k) s - x_\alpha t - x_\beta s \\
- (k + 1) r = ( (-l_\beta)(-m_\alpha) + (-k - 1)(-m_\beta - 1) ) ts - m_\alpha(-l_\alpha - 1 - k - 1) \frac{t (t - 1)}{2} \\ + t - l_\beta(-m_\beta - 1 - (-k - 1)) \frac{s (s - 1)}{2} + (m_\alpha(l_\alpha + k + 1) - x_\alpha) t + (l_\beta(m_\beta - k) - x_\beta) s - (k + 1) r.
\end{multline*}
The statement of this lemma follows from the calculations presented above.
\hfill $\square$

\subsection{Classification of Rota---Baxter Operators on $E_3$}

Throughout this entire section, we assume that $G = E_3$. Our goal is to classify the Rota---Baxter operators on $G$ up to conjugation by automorphisms $\Aut(G)$ and the action $ \ \widetilde{} \ $.

In this section, we will write $\chi = \chi_E$ for brevity.

\textbf{Lemma 6}.
Let $B$ be a Rota---Baxter operator on $G$ and let $\chi(B) = [l_\alpha, l_\beta, m_\alpha, m_\beta, x_\alpha,$ $x_\beta, 0]$. Then there exists a unique set of parameters $l_\alpha', l_\beta', m_\alpha', m_\beta', x_\alpha', x_\beta'$ satisfying the conditions
$0 \leq l_\alpha', l_\beta', x_\alpha',$ $x_\beta' \leq p-1$, 
$\, p^2 - p -1 \leq m_\alpha', m_\beta' \leq p^2$, and 
$\chi(B) = [l_\alpha', l_\beta', m_\alpha', m_\beta', x_\alpha', x_\beta', 0]$.

{\sc Proof}.
We prove existence. We choose $l_\alpha' = l_\alpha \ \mathrm{rem \, p}$, $l_\beta' = l_\beta \ \mathrm{rem \, p}$, $m_\alpha' = p^2 - l_\alpha'$, $m_\beta' = p^2 - l_\beta'$, $x_\alpha' = x_\alpha + [l_\alpha/p] - [(m_\alpha' - m_\alpha)/p] \ \mathrm{rem  \, p}$, and $x_\beta' = x_\beta + [l_\beta/p] - [(m_\beta' - m_\beta)/p] \ \mathrm{rem \, p}$. Since $a^p = b^p = c$, the new set of parameters defines the same Rota---Baxter operator $B$; thus, we obtain the required set of parameters.

We show uniqueness. Assume that $l_\alpha', l_\beta', m_\alpha', m_\beta', x_\alpha', x_\beta'$ and $l_\alpha'', l_\beta'', m_\alpha'', m_\beta'', x_\alpha'', x_\beta''$ are sets of parameters as in the statement of the theorem. Then 
$$
a^{l_\alpha'} b^{m_\alpha' \ \mathrm{rem} \ p} c^{(\ldots)} = B(a) = a^{l_\alpha''} b^{m_\alpha'' \ \mathrm{rem} \ p} c^{(\ldots)}, 
$$
consequently, $l_\alpha' = l_\alpha''$ and $m_\alpha' = m_\alpha''$. In a similar way, we obtain that $l_\beta' = l_\beta''$ and $m_\beta' = m_\beta''$. Finally, $x_\alpha' = x_\alpha'' \pmod{p}$ and $x_\beta' = x_\beta'' \pmod{p}$ by \eqref{B_gamma}; therefore, $x_\alpha' = x_\alpha''$ and $x_\beta' = x_\beta''$.
\hfill $\square$

We clarify the form of $B^{(\varphi)}$ for some $B \in RB(G)$, where $\chi(B) = [l_\alpha, l_\beta, m_\alpha, m_\beta,$ $x_\alpha, x_\beta, k]$ and $\varphi \in \Aut(G)$. 
Since $k = 0$ or $k = -1$, up to the action \ $\widetilde{}$\ \ we can assume that $k = 0$, $l_\alpha = -m_\alpha$, $l_\beta = -m_\beta$, $u_\beta = u_\alpha - 1$, $v_\alpha = \Delta - u_\alpha$, $v_\beta = \Delta - u_\alpha + 1$, and $\Delta \in \mathbb{Z}_{p^2}^{*}$. We write the formulas for the parameters of $B^{(\varphi)}$:
\begin{gather}
l_\alpha^{(\varphi)} = \frac{1}{\Delta}(u_\alpha v_\beta l_\alpha + v_\alpha v_\beta l_\beta + u_\alpha u_\beta l_\alpha + u_\beta v_\alpha l_\beta)  = u_\alpha (l_\alpha - l_\beta) + \Delta l_\beta, \nonumber \allowdisplaybreaks \\
l_\beta^{(\varphi)} = \frac{1}{\Delta}(u_\beta v_\beta l_\alpha + v_\beta^2 l_\beta + u_\beta^2 l_\alpha + u_\beta v_\beta l_\beta) = (u_\alpha - 1) (l_\alpha - l_\beta) + \Delta l_\beta, \nonumber \\ 
m_\alpha^{(\varphi)} = -\frac{1}{\Delta}(u_\alpha v_\alpha l_\alpha + v_\alpha^2 l_\beta + u_\alpha^2 l_\alpha + u_\alpha v_\alpha l_\beta) = -( u_\alpha (l_\alpha - l_\beta) + \Delta l_\beta ), \label{coeff_aut_E3} \\
m_\beta^{(\varphi)} = -\frac{1}{\Delta}(u_\beta v_\alpha l_\alpha + v_\alpha v_\beta l_\beta + u_\alpha u_\beta l_\alpha + u_\alpha v_\beta l_\beta) = -( (u_\alpha - 1) (l_\alpha - l_\beta) + \Delta l_\beta ), \nonumber \\ 
x_\alpha^{(\varphi)} = R_\alpha, \quad x_\beta^{(\varphi)} = R_\beta, \quad k^{(\varphi)} = k = 0, \nonumber
\end{gather}
where $R_\alpha$ and $R_\beta$ are defined in Appendix B by formulas \eqref{R_a} and \eqref{R_b}, respectively. One can find the details of the computations of $B^{(\varphi)}$ in Appendix B.

\textbf{Lemma 7}.
a) If $l_\alpha = l_\beta = 0$, then for any $\varphi \in \Aut(G)$, 
the equality $\chi(B^{(\varphi)}) = [0, 0, 0, 0,$ $x_\alpha^{(\varphi)}, x_\beta^{(\varphi)}, 0]$ holds.

b) If $l_\alpha = l_\beta \neq 0$, then there exists $\varphi \in \Aut(G)$ such that
$\chi(B^{(\varphi)}) = [1, 1, -1, -1,$ $x_\alpha^{(\varphi)}, x_\beta^{(\varphi)}, 0]$.

c) If $l_\alpha \neq l_\beta$, then there exists $\varphi \in \Aut(G)$ such that $\chi(B^{(\varphi)}) = [l_\alpha - l_\beta, 0, l_\beta - l_\alpha, 0, x_\alpha^{(\varphi)}$, $x_\beta^{(\varphi)}, 0]$.

{\sc Proof.}
a) The statement immediately follows from \eqref{coeff_aut_E3}.

b) We take $\Delta = l_\beta^{-1}$ in \eqref{coeff_aut_E3}.

c) We take $u_\alpha = -\Delta l_\beta (l_\alpha - l_\beta)^{-1} + 1$ in \eqref{coeff_aut_E3} for any $\Delta \neq 0$.
\hfill $\square$

\textbf{Lemma 8}.
Let $B$ be Rota---Baxter operator. Then

a) If $\chi(B) = [0, 0, 0, 0, x_\alpha, x_\beta, 0]$ and $\varphi \in \Aut(G)$, 
then $\chi(B^{(\varphi)}) = [0, 0, 0, 0, x_\alpha^{(\varphi)}, x_\beta^{(\varphi)}$, $0]$ and
$$
x_\alpha^{(\varphi)} = x_\beta + \Delta^{-1} u_\alpha (x_\alpha - x_\beta), \quad
x_\beta^{(\varphi)} = x_\beta + \Delta^{-1} (u_\alpha - 1) (x_\alpha - x_\beta).
$$

b) If $\chi(B) = [1, 1, -1, -1, x_\alpha, x_\beta, 0]$ and $\varphi \in \Aut(G)$, 
such that $\chi(B^{(\varphi)}) = [1, 1$, $-1, -1, x_\alpha^{(\varphi)}, x_\beta^{(\varphi)}, 0]$, then
\begin{gather*}
\Delta = 1, \quad
x_\alpha^{(\varphi)} = -w_\alpha + w_\beta + u_\alpha x_\alpha + (1-u_\alpha)(x_\beta + 1), \\
x_\beta^{(\varphi)} = -w_\alpha + w_\beta + (u_\alpha - 1)(x_\alpha - 1) + (2 - u_\alpha) x_\beta.
\end{gather*}

c) If $\chi(B) = [l_\alpha, 0, -l_\alpha, 0, x_\alpha, x_\beta, 0]$ 
and $\varphi \in \Aut(G)$, such that
$\chi(B^{(\varphi)}) = [l_\alpha, 0$, $-l_\alpha, 0, x_\alpha^{(\varphi)}, x_\beta^{(\varphi)}, 0]$, then 
$$
u_\alpha = 1, \quad 
x_\beta^{(\varphi)} = x_\beta, \quad
x_\alpha^{(\varphi)} = \frac{l_\alpha (1 - l_\alpha)}{2} + x_\beta + \Delta^{-1} \Big( l_\alpha(w_\beta - w_\alpha) - \left( \frac{l_\alpha(1 - l_\alpha)}{2} + x_\beta \right) + x_\alpha \Big).
$$

{\sc Proof.}
a) We compute $R_\alpha$, $R_\beta$ with $l_\alpha = l_\beta = m_\alpha = m_\beta = 0$ and, thus, obtain the required formulas for $x_\alpha^{(\varphi)}$, $x_\beta^{(\varphi)}$.

b) We obtain $\Delta = 1$ since $l_\alpha^{(\varphi)} = l_\beta^{(\varphi)} = \Delta$ by \eqref{coeff_aut_E3}. We compute $R_\alpha$, $R_\beta$ using the fact that $\Delta = 1$ and obtain
\begin{gather*}
R_\alpha = -w_\alpha + w_\beta + u_\alpha x_\alpha + (1-u_\alpha)(x_\beta + 1), \\
R_\beta = -w_\alpha + w_\beta + (u_\alpha - 1)(x_\alpha - 1) + (2 - u_\alpha) x_\beta.
\end{gather*}
Finally, we compute
$$
R_\alpha - R_\beta = u_\alpha x_\alpha + x_\beta + 1 - u_\alpha x_\beta - u_\alpha + x_\alpha - 1 - u_\alpha x_\alpha + u_\alpha + u_\alpha + u_\alpha x_\beta - 2 x_\beta = x_\alpha - x_\beta
$$

c) We have $u_\alpha = 1$ by~\eqref{coeff_aut_E3}.
We compute $R_\alpha$, $R_\beta$ under the condition that $u_\alpha = 1$ and obtain
$$
R_\alpha = \frac{l_\alpha (1 - l_\alpha)}{2} + x_\beta + \Delta^{-1} \Big( l_\alpha(w_\beta - w_\alpha) - \frac{l_\alpha(1 - l_\alpha)}{2} - x_\beta + x_\alpha \Big), \quad 
R_\beta = x_\beta .
$$
The lemma is proved.
\hfill $\square$

Now we define some families of Rota---Baxter operators on $G$:
\begin{gather*}
\mathfrak{B}^1, \quad \chi(\mathfrak{B}^1) = [0, 0, 0, 0, 1, 0, 0]; \\
\mathfrak{B}^2_x, \quad \chi(\mathfrak{B}^2_x) = [0, 0, 0, 0, x, x, 0], \quad 0 \leq x \leq p-1; \\
\mathfrak{B}^3_x, \quad \chi(\mathfrak{B}^3_x) = [1, 1, -1, -1, x, 0, 0], \quad 0 \leq x \leq p-1; \\
\mathfrak{B}^4_{l, y}, \quad \chi(\mathfrak{B}^4_{l, y}) = [l, 0, -l, 0, 0, y, 0], \quad 1 \leq l \leq p-1, \quad 0 \leq y \leq p-1.
\end{gather*}

Also the notation $\mathfrak{B}^2, \mathfrak{B}^3, \mathfrak{B}^4$ without subscripts means that we consider the specified Rota---Baxter operators with all possible parameters.

\textbf{Theorem 2}.
Any Rota---Baxter operator on $E_3$, up to conjugation by automorphisms of the group $E_3$ and the action \ $\widetilde{}$ \, is equal to one of the operators $\mathfrak{B}^1$, $\mathfrak{B}^2_x$, $\mathfrak{B}^3_x$, $\mathfrak{B}^4_{l, y}$. 
Moreover, the operators $\mathfrak{B}^1$, $\mathfrak{B}^2_x$, $\mathfrak{B}^3_x$, $\mathfrak{B}^4_{l, y}$ are contained in different orbits.

{\sc Proof}.
Let $\chi(B) = [l_\alpha, l_\beta, m_\alpha, m_\beta, x_\alpha, x_\beta, k]$. First, we prove the reducibility of the operator $B$ to any other of the selected families.

By~\eqref{RB_condition_E_3}, we have $k = 0,-1$ modulo ~$p$.
We can assume that $k = 0$, otherwise we can apply $\, \widetilde{} \,$.

{\sc Case 1}: 
Let $l_\alpha = l_\beta = 0$. If $x_\alpha = x_\beta$, then $B = \mathfrak{B}^2_{x_\alpha}$. If $x_\alpha \neq x_\beta$, then we choose $\varphi \in \Aut(G)$ with parameters $u_\alpha = 1 - x_\beta$, $\Delta = x_\alpha - x_\beta \in \mathbb{Z}_{p^2}^{*}$, and obtain $x_\alpha^{(\varphi)} = 1, x_\beta^{(\varphi)} = 0$ by Lemma 8a); consequently, $B^{(\varphi)} = \mathfrak{B}^{1}$.

{\sc Case 2}: 
Let $l_\alpha = l_\beta \neq 0$. Then there exists $\psi \in \Aut(G)$ such that $\chi(B^{(\psi)}) = [1, 1, -1, -1$, $x_\alpha^{(\psi)}, x_\beta^{(\psi)}, 0]$ by Lemma 7b). Consider $\varphi \in \Aut(G)$ with parameters $u_\alpha = 1$, $\Delta = 1$, $w_\alpha = x_\beta$, and $w_\beta = 0$; then $\chi(B^{(\varphi \psi)}) = [1, 1, -1, -1, x_\alpha^{(\psi)} - x_\beta^{(\psi)}, 0, 0]$ by Lemma~8b) and~\eqref{coeff_aut_E3}. 
Consequently, $B^{(\varphi \psi)} = \mathfrak{B}^3_{x_\alpha^{(\psi)} - x_\beta^{(\psi)} \ \mathrm{rem} \ p}$ by Lemma 6.

{\sc Case 3:} 
Let $l_\alpha \neq l_\beta$. 
Then there exists $\psi \in \Aut(G)$ such that $\chi(B^{(\psi)}) = [l_\alpha - l_\beta, 0, l_\beta - l_\alpha, 0, x_\alpha^{(\psi)}, x_\beta^{(\psi)}, 0]$ by Lemma 7b). 
We choose $\varphi \in \Aut(G)$ with parameters
\begin{gather*}
u_\alpha = 1, \quad \Delta^{-1} = -\frac{(l_\alpha - l_\beta) (1 - l_\alpha + l_\beta)}{2} - x_\beta^{(\psi)}, \\
w_\alpha = (l_\alpha - l_\beta)^{-1} \bigg( -\left( \frac{(l_\alpha - l_\beta) (1 - l_\alpha + l_\beta)}{2}  + x_\beta^{(\psi)} \right) + x_\alpha^{(\psi)} - 1 \bigg), \quad w_\beta = 0,
\end{gather*}
if $\frac{(l_\alpha - l_\beta) (1 - l_\alpha + l_\beta)}{2} + x_\beta^{(\psi)} \neq 0$, or
$$
u_\alpha = 1, \quad \Delta = 1, \quad w_\alpha = (l_\alpha - l_\beta)^{-1} x_\alpha^{(\psi)}, \quad w_\beta = 0,
$$
if $\frac{(l_\alpha - l_\beta) (1 - l_\alpha + l_\beta)}{2} + x_\beta^{(\psi)} = 0$.

Then $\chi(B^{(\varphi \psi)}) = [l_\alpha - l_\beta, 0, l_\beta - l_\alpha, 0, 0, x_\beta^{(\psi)}, 0]$ by Lemma 8c) and \eqref{coeff_aut_E3}. Consequently, $B^{(\varphi \psi)} = \mathfrak{B}^4_{l_\alpha - l_\beta \ \mathrm{rem} \ p, \ x_\beta^{(\psi)} \ \mathrm{rem} \ p}$ by Lemma 6.

Now we show uniqueness, namely, that the orbits of the operators $\mathfrak{B}^1$, $\mathfrak{B}^2_x$, $\mathfrak{B}^3_x$, $\mathfrak{B}^4_{l, y}$ under the conjugation action do not intersect.

The orbits of the operators $\mathfrak{B}^1$, $\mathfrak{B}^2_{x, x}$ do not intersect with the orbits of $\mathfrak{B}^3_{x'}$, $\mathfrak{B}^4_{l', y'}$ by item a) of Lemma 8. 
The operators $\mathfrak{B}^1$ and $\mathfrak{B}^2_{x, x}$ have non-intersecting orbits by Lemma 8a) since 
$\chi(\mathfrak{B^1})_5 - \chi(\mathfrak{B^1})_6 = 1$ and
$\chi(\mathfrak{B}^2_{x, x})_5 - \chi(\mathfrak{B}^2_{x, x})_6 = 0$.
Indeed, Lemma 8a) states that if $\chi(B) = [0, 0, 0$, $0, x_\alpha, x_\beta, 0]$ and $\varphi \in \Aut(G)$, then
$x_\alpha^{(\varphi)} - x_\beta^{\varphi} = \Delta^{-1} (x_\alpha - x_\beta)$.
The orbits of the operators $\mathfrak{B}^3_{x}$ do not intersect with the orbits of $\mathfrak{B}^4_{l, y}$ due to the fact that $\chi((\mathfrak{B}^3_{x})^{(\varphi)})_2 = \Delta \neq 0 = \chi(\mathfrak{B}^4_{l, y})_2$ for any $\varphi \in \Aut(G)$ by \eqref{coeff_aut_E3}.

We will show that the orbits of the operators from the same family do not intersect. 
We have $\chi((\mathfrak{B}^2_{x, x})^{(\varphi)})_5 = x$ for any $\varphi \in \Aut(G)$ by Lemma 8a). Consequently, the orbits of $\mathfrak{B}^2_{x, x}$ and $\mathfrak{B}^2_{x', x'}$ do not intersect for $x \neq x'$. 
The orbits of $\mathfrak{B}^3_{x}$ and $\mathfrak{B}^3_{x'}$ do not intersect for $x \neq x'$ by virtue of Lemma 8b). 
Indeed, Lemma 8b) states that if $\chi(B) = [1, 1, -1, -1, x_\alpha, x_\beta, 0]$ and $\varphi \in \Aut(G)$ is such that $\chi(B^{(\varphi)}) = [1, 1$, $-1, -1, x_\alpha^{(\varphi)}, x_\beta^{(\varphi)}, 0]$, then
$x_\alpha^{(\varphi)} - x_\beta^{(\varphi)} = x_\alpha - x_\beta$.
Therefore, $\mathfrak{B}^3_{x}$ and $\mathfrak{B}^3_{x'}$ lie in different orbits since
$\chi(\mathfrak{B}^3_{x})_5 - \chi(\mathfrak{B}^3_{x})_6 = x \neq x' = \chi(\mathfrak{B}^3_{x'})_5 - \chi(\mathfrak{B}^3_{x'})_6$. 
Finally, the orbits of $\mathfrak{B}^4_{l, y}$ and $\mathfrak{B}^4_{l', y'}$ do not intersect for $l \neq l'$ or $y \neq y'$. Let $l = l'$ and $y \neq y'$; then we obtain the required result by Lemma 8c). Let $l \neq l'$ and let $\varphi \in \Aut(G)$ preserve the second and fourth coordinates of $\chi(B^{(\varphi)})$ equal to zero. Then $\chi((\mathfrak{B}^4_{l, y})^{(\varphi)})_1 = l \neq l' = \chi((\mathfrak{B}^4_{l', y'})^{(\varphi)})_1$ by~\eqref{coeff_aut_E3}.
\hfill $\square$

\subsection{Classification of Rota---Baxter Operators on the Group $H_3$}

Throughout this section, we assume that $G = H_3$. Our goal is to classify all Rota---Baxter operators on $G$ up to conjugation by automorphisms of the group $G$ and the action $ \ \widetilde{} \ $.

By Proposition 2b), conjugating a Rota---Baxter operator by automorphisms yields a Rota---Baxter operator again. 
We denote $B^{(\varphi)} = \varphi^{-1} \circ B \circ \varphi$.

We clarify the form of the operator $B^{(\varphi)}$ for a Rota---Baxter operator $B$ such that $\chi(B) = (l_\alpha, l_\beta, m_\alpha, m_\beta, x_\alpha, x_\beta, k)$ and $\varphi \in \Aut(G)$:
\begin{equation}\label{coeff_aut_H3}
\begin{gathered}
l_\alpha^{(\varphi)} = \frac{1}{\Delta}(u_\alpha v_\beta l_\alpha + v_\alpha v_\beta l_\beta - u_\alpha u_\beta m_\alpha - u_\beta v_\alpha m_\beta), \\
l_\beta^{(\varphi)} = \frac{1}{\Delta}(u_\beta v_\beta l_\alpha + v_\beta^2 l_\beta - u_\beta^2 m_\alpha - u_\beta v_\beta m_\beta), \\ 
m_\alpha^{(\varphi)} = \frac{1}{\Delta}(-u_\alpha v_\alpha l_\alpha - v_\alpha^2 l_\beta + u_\alpha^2 m_\alpha + u_\alpha v_\alpha m_\beta), \\
m_\beta^{(\varphi)} = \frac{1}{\Delta}(-u_\beta v_\alpha l_\alpha - v_\alpha v_\beta l_\beta + u_\alpha u_\beta m_\alpha + u_\alpha v_\beta m_\beta), \\ 
x_\alpha^{(\varphi)} = R_\alpha, \quad x_\beta^{(\varphi)} = R_\beta, \quad k^{(\varphi)} = k,
\end{gathered}
\end{equation}
where $R_\alpha$ and $R_\beta$ are defined in Appendix B by formulas \eqref{R_a} and \eqref{R_b}, respectively. One can find the details of the computation of $B^{(\varphi)}$ in Appendix B.

We note some connections between the parameters of the operators $B^{(\varphi)}$ and $B$.

\textbf{Lemma 9}.
Let $\chi(B) = (l_\alpha, l_\beta, m_\alpha, m_\beta, x_\alpha, x_\beta, k)$. Then

a) $l_\alpha^{(\varphi)} + m_\beta^{(\varphi)} = l_\alpha + m_\beta$ and $k$ = $k^{(\varphi)}$ for any $\varphi \in \Aut(G)$;

b) an automorphism $\varphi$ such that $l_\alpha^{(\varphi)} = 0$ exists, $m_\beta^{(\varphi)} = l_\alpha + m_\beta$, if at least one of the following conditions holds:
1) $(l_\beta, m_\alpha)\neq(0,0)$,
2) $l_\alpha \neq m_\beta$;

c) let $l_\alpha = 0$ and $l_\alpha^{(\varphi)} = 0$
for any $\varphi \in \Aut(G)$. Then 
$l_\beta^{(\varphi)} m_\alpha^{(\varphi)} = l_\beta m_\alpha = -k (m_\beta + 1)$;

d) if $l_\alpha = m_\beta$ and $l_\beta = m_\alpha = 0$, then $l_\alpha^{(\varphi)} = m_\beta^{(\varphi)} = l_\alpha$ for any $\varphi \in \Aut(G)$.

{\sc Proof.}
a) From \eqref{coeff_aut_H3}, we obtain that $l_\alpha^{(\varphi)} + m_\beta^{(\varphi)} = l_\alpha + m_\beta$.

b)
We consider three cases. Let $m_\alpha \neq 0$; then we choose $u_\alpha = v_\beta = 1$ and $v_\alpha = 0$, and we obtain $l_\alpha^{(\varphi)} = l_\alpha - t_\beta m_\alpha$. Now we choose $t_\beta = l_\alpha m_\alpha^{-1}$, and we have $l_\alpha^{(\varphi)} = 0$. Let us consider the case $l_\beta \neq 0$. Let $t_\alpha = s_\beta = 1$ and $t_\beta = 0$; then $l_\alpha^{(\varphi)} = l_\alpha + s_\alpha l_\beta$. Now we choose $s_\alpha = -l_\alpha l_\beta^{-1}$, then $l_\alpha^{(\varphi)} = 0$.

We analyze the last case. Let $l_\beta = m_\alpha = 0$, $l_\alpha \neq m_\beta$, and $u_\alpha = u_\beta = 1$; then $l_\alpha^{(\varphi)} = \frac{1}{v_\beta - v_\alpha}(v_\beta l_\alpha - v_\alpha m_\beta)$. Now we consider $v_\alpha = l_\alpha$ and $v_\beta = m_\beta$; we can do this since $l_\alpha \neq m_\beta$, and we obtain $l_\alpha^{(\varphi)} = 0$.

c) By~\eqref{RB_condition} and~\eqref{coeff_aut_H3}, we have 
$-l_\beta m_\alpha = k (m_\beta + 1)$ and $-l_\beta^{(\varphi)} m_\alpha^{(\varphi)} = k (l_\alpha^{(\varphi)} + m_\beta^{(\varphi)} + 1) = k (m_\beta^{(\varphi)} + 1)$.
Thus, we obtain $l_\beta m_\alpha = l_\beta^{(\varphi)} m_\alpha^{(\varphi)} = -k (m_\beta + 1)$ by a).

d) This point follows immediately from~\eqref{coeff_aut_H3}.
\hfill $\square$

\textbf{Lemma 10}.
Let $\chi(B) = (l_\alpha, l_\beta, m_\alpha, m_\beta, x_\alpha, x_\beta, k)$. Then for $k \neq 0, -1$ there exist $\varphi \in \Aut(G)$ such that $\chi(B^{(\varphi)}) = (l_\alpha, l_\beta, m_\alpha, m_\beta, 0, 0, k)$.

{\sc Proof.}
Let $\varphi \in \Aut(G)$ with parameters $t_\beta = s_\alpha = 0$ and $t_\alpha = s_\beta = 1$.
Then by~\eqref{coeff_aut_H3}, we have
\begin{gather*}
\Delta = 1, \quad l_\alpha^{(\varphi)} = l_\alpha, \quad l_\beta^{(\varphi)} = l_\beta, \quad m_\alpha^{(\varphi)} = m_\alpha, \quad m_\beta^{(\varphi)} = m_\beta, \\
x_\alpha^{(\varphi)} = (k - l_\alpha) w_\alpha - m_\alpha w_\beta + x_\alpha, \quad
x_\beta^{(\varphi)} =  -l_\beta w_\alpha + (k - m_\beta) w_\beta + x_\beta.
\end{gather*}
We compute the determinant of the system of linear equations $x_\alpha^{(\varphi)} = 0$, $x_\beta^{(\varphi)} = 0$ with respect to the variables $w_\alpha, w_\beta$:
\begin{gather*}
\begin{vmatrix}
(k - l_\alpha) & -m_\alpha \\
-l_\beta & (k - m_\beta)
\end{vmatrix}
 = k^2 + l_\alpha m_\beta - k l_\alpha - l_\beta m_\alpha - km_\beta
 \stackrel{\eqref{RB_condition}}{=} k^2 + k.
\end{gather*}
The determinant is equal to zero if and only if $k = 0, -1$.
That is, in the case when $k \neq 0, -1$, we can choose $w_\alpha, w_\beta$ such that $x_\alpha^{(\varphi)} = 0$, $x_\beta^{(\varphi)} = 0$.
\hfill $\square$

\textbf{Proposition 6}.
Let $B$ be a Rota---Baxter operator, $\chi(B) = (0, l_\beta, m_\alpha$, $m_\beta, x_\alpha,$ $x_\beta, k)$, and $(l_\beta, m_\alpha, m_\beta) \neq (0, 0, 0)$. 

a) If $k \neq 0$ and $m_\beta \neq -1$, then there exists $\varphi \in \Aut(G)$ such that $\chi(B^{(\varphi)}) = (0, -1, k(m_\beta + 1), m_\beta, x_\alpha^{(\varphi)}, x_\beta^{(\varphi)}, k)$.

b) If $k = m_\beta = 0$, then there exists $\varphi \in \Aut(G)$ such that $\chi(B^{(\varphi)}) = (0, 1, 0, 0, x_\alpha^{(\varphi)}, \linebreak x_\beta^{(\varphi)}, 0)$.

c) If $k = 0$ and $m_\beta \neq 0$, then there exists $\varphi \in \Aut(G)$ such that $\chi(B^{(\varphi)}) = (0, 0, 0, m_\beta$, $x_\alpha^{(\varphi)},x_\beta^{(\varphi)}, 0)$.

d) If $m_\beta = -1$, then there exists $\varphi \in \Aut(G)$ such that $\chi(B^{(\varphi)}) = (0, 0, 0, -1, x_\alpha^{(\varphi)}, \linebreak x_\beta^{(\varphi)}, k)$.

{\sc Proof}.
Below in the proof, we consider automorphisms~$\varphi$ such that $l_\alpha^{(\varphi)} = 0$. From item c) of Lemma 9, we have $l_\beta^{(\varphi)} m_\alpha^{(\varphi)} = l_\beta m_\alpha = -k (m_\beta + 1)$. 
We use formula \eqref{coeff_aut_H3} to prove all items of this proposition.

a) Let $l_\beta \neq 0$. We consider an automorphism $\varphi$ with parameters $u_\alpha = l_\beta$, $u_\beta = 0$, $v_\alpha = 0$, $v_\beta = -1$, and arbitrary $w_\alpha$, $w_\beta$.

Now let $l_\beta = 0$ and $m_\alpha \neq 0$. We consider an automorphism $\varphi$ with parameters $u_\alpha = -m_\beta$, $u_\beta = -1$, $v_\alpha = m_\alpha$, $v_\beta = 0$, and arbitrary $w_\alpha$, $w_\beta$. Then $B^{(\varphi)}$ has the required form.

From the conditions $k \neq 0$ and $m_\beta \neq -1$, we obtain that $l_\beta$ and $m_\alpha$ are not equal to zero simultaneously. Thus, item a) is proved.

b) We have $l_\beta m_\alpha = 0$ and $(l_\beta, m_\alpha) \neq (0, 0)$ by assumption; consequently, exactly one of the parameters $l_\beta, m_\alpha$ is equal to zero.

If $l_\beta = 0$ and $m_\alpha \neq 0$, we choose $\varphi$ with parameters $u_\alpha = 0$, $u_\beta = m_\alpha^{-1}$, $v_\alpha = 1$, $v_\beta = 0$, and arbitrary~$w_\alpha$, $w_\beta$.

If $l_\beta \neq 0$ and $m_\alpha = 0$, we consider an automorphism $\varphi$ with parameters $u_\alpha = l_\beta$, $u_\beta = 0$, $v_\alpha = 0$, $v_\beta = 1$, and arbitrary $w_\alpha$, $w_\beta$.

c) We have $l_\beta m_\alpha = 0$; therefore, at least one of the parameters $l_\beta, m_\alpha$ is equal to zero.

Let $l_\beta = 0$; in this case, we consider an automorphism $\varphi$ with parameters $u_\alpha = -m_\beta$, $u_\beta = 0$, $v_\alpha = m_\alpha$, $v_\beta = 1$, and arbitrary $w_\alpha$, $w_\beta$.

Now let $m_\alpha = 0$; we choose an automorphism $\varphi$ with parameters $u_\alpha = 1$, $u_\beta = l_\beta$, $v_\alpha = 0$, $v_\beta = m_\beta$, and arbitrary $w_\alpha$, $w_\beta$.

d) In this case, the condition $l_\beta m_\alpha = 0$ holds again; we choose the automorphism from the proof of item c) and obtain the required result.
\hfill $\square$

We define some families of RB-operators:
\begin{gather*}
D^1_{x, y}, \quad \chi(D^1_{x, y}) = (0, 0, 0, 0, x, y, 0); \\
D^2_{x, y}, \quad \chi(D^2_{x, y}) = (0, 1, 0, 0, x, y, 0); \\
D^3_{m, x, y}, \quad \chi(D^3_{m, x, y}) = (0, 0, 0, m, x, y, 0), \quad m \neq 0; \allowdisplaybreaks \\
B^1_{k, m}, \quad \chi(B^1_{k, m}) = (0, -1, k(m + 1), m, 0, 0, k), \quad k \neq 0, -1, \quad m \neq -1; \\
B^2_k, \quad \chi(B^2_k) = (0, 0, 0, -1, 0, 0, k), \quad k \neq 0, -1. \\
B^3_m, \quad \chi(B^3_m) = (m, 0, 0, m, 0, 0, m^2/(2m+1)), \quad m \neq 0, -1/2. \end{gather*}

Here, the notation $D^1, D^2, D^3, B^1, B^2, B^3$ without subscripts means that we consider the specified Rota---Baxter operators with all possible parameters.

\textbf{Proposition 7}.
The families $D^1, D^2, D^3, B^1, B^2, B^3$ have non-intersecting orbits under the conjugation action by automorphisms of $G$.

{\sc Proof}. 
The orbits of the families $D^1$--$D^3$ and the orbits of $B^1$--$B^3$ do not intersect by Lemma 9a) and because $\chi(D^i) \neq \chi(B^j)$ for $i, j \in \{1, 2, 3\}$.

The orbits of the family $D^1$ do not intersect with the orbits of the families $D^2, D^3$ since \linebreak
$\chi((D^1)^{(\varphi)}) = (0, 0, 0, 0, x^{(\varphi)}, y^{(\varphi)}, 0)$ for any $\varphi \in \Aut(G)$ by formula \eqref{coeff_aut_H3}. The orbits of the family $D^2$ do not intersect with the orbits of the family $D^3$ by item a) of Lemma 9 and due to the fact that $\chi(D^2)_4 = 0 \neq \chi(D^3)_4$.

The orbits of the families $B^3$ do not intersect with the orbits of the families $B^1, B^2$ by item d) of Lemma~9. 
Finally, the orbits of the family $B^1$ and the orbits of the family $B^2$ do not intersect by item~a) of Lemma 9 and because $\chi(B^1)_4 \neq -1 = \chi(B^2)_4$.
\hfill $\square$

We have already shown that the orbits of the families $B^1$, $B^2$, $B^3$ do not intersect with each other. Now we prove that the orbits of the operators within each family do not intersect.

\textbf{Proposition 8}.

a) For any $k_1 \neq 0, -1$, $m_1 \neq -1$ and any $k_2 \neq 0, -1$, $m_2 \neq -1$ such that $k_1 \neq k_2$ or $m_1 \neq m_2$, there does not exist $\varphi \in \Aut(G)$ such that $(B^1_{m_1, k_1})^{(\varphi)} = B^1_{m_2, k_2}$.

b) For any $k_1 \neq 0, -1$ and $k_2 \neq 0, -1$ such that $k_1 \neq k_2$, there does not exist an automorphism $\varphi \in \Aut(G)$ such that $(B^2_{k_1})^{(\varphi)} = B^2_{k_2}$.

c) For any $m_1$ and $m_2$ such that $m_1 \neq m_2$ and $m_1, m_2 \neq 0, -1/2$, there does not exist $\varphi \in \Aut(G)$ such that $(B^3_{m_1})^{(\varphi)} = B^3_{m_2}$.

{\sc Proof.}
All items of this proposition are a direct consequence of Lemma 9.
Item a) of this proposition is true by virtue of the statements of items a) and c) of Lemma~9. 
The proof of item b) is obvious modulo item a) of Lemma 9 and, finally, item c) of this proposition follows from item d) of Lemma 9.
\hfill $\square$

By Proposition 6 and Lemma 10, we obtain that all operators with $k \neq 0, -1$ reduce to a unique operator from the families $B^2, B^4, B^5$ by conjugation by automorphisms and the action $ \ \widetilde{} \ $, we will present this statement at the end of the section.

Now we concentrate to the families $D^1, D^2, D^3$.

\textbf{Lemma 11}.

a) An automorphism $\varphi$ maps by conjugation a Rota---Baxter operator $D^2_{x, y}$ with some $x, y$ to an operator of the same family, that is, $(D^2_{x, y})^{(\varphi)} = D^2_{x^{(\varphi)}, y^{(\varphi)}}$, if and only if $t_\alpha = s_\beta \neq 0$ and $s_\alpha = 0$.

b) An automorphism $\varphi$ maps by conjugation a Rota---Baxter operator $D^3_{x, y}$ with some $x, y$ and $m \neq 0$ to an operator of the same family, that is, $(D^3_{m, x, y})^{(\varphi)} = D^3_{m, x^{(\varphi)}, y^{(\varphi)}}$, if and only if $t_\beta = s_\alpha = 0$.

{\sc Proof}.
Sufficiency is verified by direct substitution into~\eqref{coeff_aut_H3}.
We write out the necessary conditions on the parameters $t_\alpha, t_\beta, s_\alpha, s_\beta$ using formulas~\eqref{coeff_aut_H3}.

a) We require that the conditions $l_\alpha^{(\varphi)} = m_\alpha^{(\varphi)} = m_\beta^{(\varphi)} = 0$ and $l_\beta^{(\varphi)} = 1$ hold for the desired automorphism $\varphi$. This is equivalent to the restrictions $s_\alpha s_\beta = 0$, $s_\alpha = 0$, and $t_\alpha^{-1} s_\beta = 1$ on the parameters of the automorphism $\varphi$; hence, we obtain that $t_\alpha = s_\beta \neq 0$ and $s_\alpha = 0$.

b) Necessarily, the conditions $l_\alpha^{(\varphi)} = l_\beta^{(\varphi)} = m_\alpha^{(\varphi)} = 0$ and $m_\beta^{(\varphi)} = m$ must hold for the desired automorphism $\varphi$. This is equivalent to the restrictions $t_\beta s_\alpha = 0$, $t_\beta s_\beta = 0$, and $t_\alpha s_\alpha = 0$ on the parameters of the automorphism $\varphi$; hence, we obtain that $t_\beta = 0$ and $s_\alpha = 0$ since $t_\alpha s_\beta$ must be non-zero.
\hfill $\square$

Now, conjugating by automorphisms, we reduce the operators of the families $D^1, D^2, D^3$ to a simpler form.

\textbf{Proposition 9}.

a) If at least one of the parameters $x, y$ is non-zero, then there exists $\varphi \in \Aut(G)$ such that $\chi((D^1_{x, y})^{(\varphi)}) = (0, 0, 0, 0, 1, 0, 0)$.

b) If $x \neq 0,$, then there exists $\varphi \in \Aut(G)$ such that $\chi((D^2_{x, y})^{(\varphi)}) = (0, 1, 0, 0, 1, 0, 0)$. If $x = 0$, then there exists $\varphi \in \Aut(G)$ such that $\chi((D^2_{x, y})^{(\varphi)}) = (0, 1, 0, 0, 0, 0, 0)$.

c) If $x \neq 0$, then there exists $\varphi \in \Aut(G)$ such that $\chi((D^3_{m, x, y})^{(\varphi)}) = (0, 0, 0, m, 1, 0, 0)$. If $x = 0$, then there exists $\varphi \in \Aut(G)$ such that $\chi((D^3_{m, x, y})^{(\varphi)}) = (0, 0, 0, m, 0, 0, 0)$.

{\sc Proof.}
a) From~\eqref{coeff_aut_H3}, we have
$$
x_\alpha^{(\varphi)} = \frac{1}{u_\alpha v_\beta - u_\beta v_\alpha}(u_\alpha x + v_\alpha y), \quad x_\beta^{(\varphi)} = \frac{1}{u_\alpha v_\beta - u_\beta v_\alpha}(u_\beta x + v_\beta y).
$$
If $x \neq 0$, then we take $\varphi \in \Aut(G)$ with parameters $u_\alpha = 1$, $u_\beta = -y$, $v_\alpha = 0$, and $v_\beta = x$.
If $y \neq 0$, then we take $\varphi \in \Aut(G)$ with parameters $u_\alpha = 0$, $u_\beta = -y$, $v_\alpha = 1$, and $v_\beta = x$.

b) By item a) of Lemma 11, $u_\alpha = v_\beta \neq 0$ and $v_\alpha = 0$, hence we obtain that $x_\alpha^{(\varphi)} = u_\alpha^{-1} x$ and $x_\beta^{(\varphi)} = -w_\alpha + \frac{u_\beta x + u_\alpha y}{\Delta}$.
If $x \neq 0$, then we consider an automorphism $\varphi$ with parameters $u_\alpha = x$, $u_\beta = 0$, and $w_\alpha = y/x$. If $x = 0$, then we consider an automorphism $\varphi$ with parameters $u_\alpha = 1$, $u_\beta = 0$, and $w_\alpha = y$.

c) By item b) of Lemma 11, $u_\beta = v_\alpha = 0$, hence we obtain that 
$x_\alpha^{(\varphi)} = x/v_\beta$ and 
$x_\beta^{(\varphi)} = (-w_\beta m + v_\beta y)/(u_\alpha v_\beta)$.
Let $x \neq 0$; then we consider an automorphism $\varphi$ with parameters $u_\alpha = 1$, $v_\beta = x$, and $w_\beta = xy/m$. If $x = 0$, then we consider an automorphism $\varphi$ with parameters $u_\alpha = v_\beta = 1$ and $w_\beta = y/m$.
\hfill $\square$

Let us define new families of operators, which the families can be reduced by means of the previous proposition
 $D^1, D^2, D^3$:
\begin{gather*}
B^0 = (0, 0, 0, 0, 0, 0, 0); \quad B^4 = (0, 0, 0, 0, 1, 0, 0); \\
B^5 = (0, 1, 0, 0, 1, 0, 0); \quad B^6 = (0, 1, 0, 0, 0, 0, 0); \\
B^7_m = (0, 0, 0, m, 1, 0, 0), \quad m\neq 0; \quad B^8_m = (0, 0, 0, m, 0, 0, 0), \quad m\neq 0.
\end{gather*}

\textbf{Lemma 12}.
The orbits of the Rota---Baxter operator families $B^0$--$B^8$ under the conjugation action by automorphisms do not intersect each other, and for any $m_1 \neq m_2$, there does not exist an automorphism $\varphi \in \Aut(G)$ such that $(B^i_{m_1})^{(\varphi)} = B^i_{m_2}$ for $i = 7, 8$.

{\sc Proof}.
The orbit of $B^0$ does not intersect with the orbit of $B^4$ since $(B^0)^{(\varphi)} = B^0$ for any $\varphi \in \Aut(G)$. The orbit of $B^5$ does not intersect with the orbit of $B^6$ and the orbit of $B^7$ does not intersect with the orbit of $B^8$ due to considerations of the order of the images of these operators, namely due to the fact that $|\Imm(B^5)| = p^2$, $|\Imm(B^4)| = |\Imm(B^6)| = p$, and $|\Imm(B^7_{m_1})| = p^2$, $|\Imm(B^8_{m_2})| = p$ for any $m_1, m_2 \neq 0$.

The statement that $(B^i_{m_1})^{(\varphi)} \neq B^i_{m_2}$ for any $m_1 \neq m_2$, any automorphism $\varphi \in \Aut(G)$, and $i = 7, 8$ follows from item a) of Lemma 9.

Taking into account Proposition 6, the lemma is completely proved.
\hfill $\square$

Below we formulate a proposition that establishes the reducibility to the families $B^0$--$B^8$.

\textbf{Proposition 10}.
Let $B$ be a Rota---Baxter operator, $\chi(B) = (l_\alpha, l_\beta, m_\alpha, m_\beta, x_\alpha$, $x_\beta, k)$. Then, if $k = 0$, there exists $\varphi \in \Aut(G)$ such that $B^{(\varphi)} \in B^0 \cup \underset{j = 4}{\overset{8}{\bigcup}} B^j$.

{\sc Proof}.
Let $k = 0$; then by points b) and c) of Proposition~6, there exists $\varphi \in \Aut(G)$ such that $B^{(\varphi)} \in D^1 \cup D^2 \cup D^3$. By Proposition 9, there exists $\psi \in \Aut(G)$ such that $B^{(\psi \varphi)} \in \underset{j = 4}{\overset{8}{\bigcup}} B^j$.
\hfill $\square$

\textbf{Theorem 3}.
Any non-trivial Rota---Baxter operator on $H_3$, up to conjugation by automorphisms of the group $H_3$ and the action \ $\widetilde{}$ \, is equal to one of the following operators:
\begin{gather*}
B^1_{k, m}, \quad \chi(B^1_{k, m}) = (0, -1, k(m + 1), m, 0, 0, k), \quad 1 \leq k < (p-1)/2, \quad m \in \mathbb{Z}_p\setminus \{-1\}; \\
B^2_k, \quad \chi(B^2_k) = (0, 0, 0, -1, 0, 0, k), \quad 1 \leq k < (p-1)/2; \\
B^3_m, \quad \chi(B^3_m) = (m, 0, 0, m, 0, 0, m^2/(2m+1)), \quad 1 \leq m < (p-1)/2; \allowdisplaybreaks \\
O^1_m, \quad \chi(O^1_m) = (0, -1, (p-1)/2\cdot(m + 1), m, 0, 0, (p-1)/2), \quad 0 \leq m < (p-1)/2; \allowdisplaybreaks \\
O^2, \quad \chi(O^2) = (0, 0, 0, -1, 0, 0, (p-1)/2); \allowdisplaybreaks \\
B^4 = (0, 0, 0, 0, 1, 0, 0); \quad 
B^5 = (0, 1, 0, 0, 1, 0, 0); \quad 
B^6 = (0, 1, 0, 0, 0, 0, 0); \\
B^7_m = (0, 0, 0, m, 1, 0, 0), \quad m \in \mathbb{Z}_p^*; 
\quad B^8_m = (0, 0, 0, m, 0, 0, 0), \quad m \in \mathbb{Z}_p^*.
\end{gather*}
Moreover, the operators $B^1$--$B^8$, $O^1,O^2$ are contained in different orbits.

{\sc Proof}.
Let $S\in \underset{t = 1}{\overset{8}{\bigcup}} B^t \cup O^1\cup O^2$. We show that either $(\widetilde{S})^{(\varphi)} \not \in \underset{t = 1}{\overset{8}{\bigcup}} B^t \cup O^1 \cup O^2$ for any $\varphi \in \Aut(G)$, or $(\widetilde{S})^{(\psi)} = S$ for some $\psi\in\Aut(G)$. 

Let $S\in \underset{t = 4}{\overset{8}{\bigcup}} B^t$. 
By Lemma 5, the equality $\chi(\widetilde{S})_7 = -1$ holds, 
and by item a) of Lemma 9, it is true for any $\varphi \in \Aut(G)$ that $\chi((\widetilde{S})^{(\varphi)})_7 = -1$.

One can easily verify that $(\widetilde{S})^{(\varphi)}\notin \underset{t = 1}{\overset{8}{\bigcup}} B^t \cup O^1 \cup O^2$ holds for any $\varphi \in \Aut(G)$.

Let $S\in B^1\cup B^2$. 
By Lemma 5, we compute $\frac{p-1}{2} < \chi(\widetilde{S})_7\leq p - 2$; hence, we obtain that $(\widetilde{S})^{(\varphi)}\notin \underset{t = 1}{\overset{8}{\bigcup}} B^t \cup O^1 \cup O^2$ by item a) of Lemma 9.

Let $S = B^3_m$. 
Again, by Lemma 5, we deduce $\frac{p-1}{2} <\chi(\widetilde{S})_1 \leq p - 2$ and $\chi(\widetilde{S})_1 = \chi(\widetilde{S})_4$; from item d) of Lemma 9, it follows that $(\widetilde{S})^{(\varphi)} \notin \underset{t = 1}{\overset{8}{\bigcup}} B^t \cup O^1 \cup O^2$.

Let $S = O^1_m$. Then from the previous items of the proof, we obtain that $(\widetilde{S})^{(\varphi)} \notin \underset{t = 1}{\overset{8}{\bigcup}} B^t$.
By Lemma 5, we obtain the following conditions:
$\chi(\widetilde{S})_7 = \frac{p-1}{2}$, 
$\frac{p+1}{2} \leq \chi(\widetilde{S})_4 \leq p - 1$, and 
$\chi(\widetilde{S})_1 = -1$. 
Then, by virtue of item a) of Lemma 9, there does not exist an automorphism $\varphi$ such that $\widetilde{S}^{(\varphi)} \in O^1 \cup O^2$.
We find that $(\widetilde{S})^{(\varphi)} \notin \underset{t = 1}{\overset{8}{\bigcup}} B^t \cup O^1 \cup O^2$ holds for any $\varphi\in\Aut(G)$.

Let $S = O^2$. 
By Lemma 5, the relations $\chi(\widetilde{S})_7 = \frac{p-1}{2}$, $\chi(\widetilde{S})_1 = -1$, and $\chi(\widetilde{S})_4 = 0$ hold. 
Then, by Lemma 9 b) and Proposition 6 d), there exists $\varphi \in \Aut(G)$ such that $(\widetilde{S})^{(\varphi)} = O^2$.

By Proposition 8 and Lemma 12, we find that no operator from $\underset{t = 1}{\overset{8}{\bigcup}} B^t \cup O^1 \cup O^2$ is conjugate to any other operator from this union.

Now we show that an arbitrary non-trivial RB-operator $B$, $\chi(B) = (l_\alpha, l_\beta$, $m_\alpha,m_\beta, x_\alpha$, $x_\beta, k)$, reduces to an operator from $\underset{t = 1}{\overset{8}{\bigcup}} B^t \cup O^1 \cup O^2$ by conjugation by automorphisms and the application of \, $\widetilde{}$\,.

{\sc Case 1}: $k = 0$. 
By Proposition 10, there exists $\varphi \in \Aut(G)$ such that $B^{(\varphi)} \in \underset{t = 1}{\overset{8}{\bigcup}} B^t \cup O^1 \cup O^2$.

{\sc Case 2}: $k = -1$. Applying \ $\widetilde{}$\, we obtain by Lemma 5 that $\chi(\widetilde{B})_7 = 0$, which is Case~1.

{\sc Case 3}: $k \neq 0, -1$. 

{\sc Case 3a}: $l_\alpha \neq m_\beta$ or $(l_\beta, m_\alpha)\neq(0,0)$. 
If $k < \frac{p-1}{2}$, then by item b) of Lemma 9, Proposition 6, and Lemma 10, there exists $\varphi \in \Aut(G)$ such that $B^{(\varphi)} \in B^1 \cup B^2$.
If $k > \frac{p-1}{2}$, then after applying \ $\widetilde{}$\ \ to $B$, we obtain by Lemma~5 that $1 \leq \chi(\widetilde{B})_7 < \frac{p-1}{2}$, and this case reduces to the previous one.
Now let $k = \frac{p-1}{2}$. If $0 \leq l_\alpha + m_\beta < \frac{p-1}{2}$, then by item b) of Lemma 9, Proposition~5, and~Lemma 10, there exists $\varphi \in \Aut(G)$ such that $B^{(\varphi)} \in O^1$. If $\frac{p-1}{2} \leq l_\alpha + m_\beta < p-1$, then $0 \leq \chi(\widetilde{B})_1 + \chi(\widetilde{B})_4 < \frac{p-1}{2}$, and this case reduces to the previous one. If $l_\alpha + m_\beta = p-1$, then by item b) of Lemma 9, Proposition~5, and~Lemma 10, there exists $\varphi \in \Aut(G)$ such that $B^{(\varphi)} \in O^2$.

{\sc Case 3b}: $l_\alpha = m_\beta$ and $l_\beta = m_\alpha = 0$. 
Then $k = \frac{m_\beta^2}{2 m_\beta + 1}$ by~\eqref{RB_condition}, and at the same time $m_\beta \neq 0, \frac{p-1}{2}$. 
Up to the action \ $\widetilde{}$\, we can assume that $1 \leq m_\beta < \frac{p-1}{2}$.
In this case, by Lemma~9d) and~Lemma~10, there exists $\varphi \in \Aut(G)$ such that $B^{(\varphi)} = B^3_{m_\beta}$.
\hfill $\square$

\section{Examples and Properties of RB-Operators on an Extraspecial Group}

\subsection{Analogue of the Spectral Property}

Since Rota---Baxter operators on groups are analogues of Rota---Baxter operators on algebras, it is logical to ask which properties can be transferred and in which cases. For Rota---Baxter operators, the following result holds.

\textbf{Theorem 1}~\cite{G}.
Let $A$ be a unital finite-dimensional algebra over a field $\mathbb{F}$ and let $\lambda\in F$. 
Then for any Rota---Baxter operator $R$ of weight $\lambda$ on $A$, we have
$\Spec(R)\subset \{0,-\lambda\}$.

This theorem implies the existence of an invariant $\mathrm{rb}_{\lambda}(A)$ on a unital finite-dimensional algebra $A$, which is the minimal $N \in \mathbb{N}$ such that for any Rota---Baxter operator $R$ of weight $\lambda$ on $A$, there exists $m \leq N$ for which $R^{m} (R + \lambda \ \text{id})^{N-m} = 0$ holds.

\textbf{Example}~\cite{G}.
Let $\lambda \in \mathbb{C}$. Then $rb_\lambda(M_n(\mathbb{C})) = 2n-1$.

For which finite groups will the analogue of this property fail to hold? 
Equivalently, for which finite groups does there exist a Rota---Baxter operator $B$ for any $N \in \mathbb{N}$ such that $B^m \widetilde{B}^{N-m} \neq \{ e \}$ for all $m \leq N$?

If the analogue of the spectral property does not hold, this means that for some $d_1, d_2 \in \mathbb{N}$ and a Rota---Baxter operator $B$, the operators $B$ and $\widetilde{B}$ are bijective on $\Imm(B^{d_1} \widetilde{B}^{d_2}) \neq \{ e \}$. Then it is reasonable to ask for which finite groups $G$ there exists a Rota---Baxter operator such that $B$ and $\widetilde{B}$ are bijective. For non-abelian groups, the smallest example by order is the group $H_3$ for $p = 3$.

We describe all Rota---Baxter operators $B$ on $H_3$ such that $B$ and $\widetilde{B}$ are bijections, and~we provide a construction for these operators.

\subsection{Bijective Operators on the Heisenberg Group}

\textbf{Proposition 11}.
Let $B$ be a Rota---Baxter operator on $H_3$. Then $B$ and $\widetilde{B}$ are bijective if and only if the operator $B$, up to conjugation by automorphisms and the action \ $\widetilde{}$\, is contained in $B^1 \cup B^3 \cup O^1$.

{\sc Proof.}
Let $\chi(B) = (l_\alpha, l_\beta, m_\alpha, m_\beta, x_\alpha, x_\beta, k)$. One can easily understand from condition \eqref{RB_condition} that the operator $B$ is surjective if and only if $l_\alpha + m_\beta \neq -1$ and $k \neq 0$; moreover, the surjectivity of the mapping is equivalent to bijectivity since $G$ is finite.

Assume that the operators $B$ and $\widetilde{B}$ are bijective; then $k \neq 0, -1$ by item a) of Lemma 9 and Lemma 5. Then there exists $\varphi \in \Aut(G)$ such that $B^{(\varphi)} \in B^1 \cup B^3 \cup O^1$ up to the action \ $\widetilde{}$\ , moreover, $B^{(\varphi)}$ cannot be contained in the families $B^2$ and $O^2$ since $\chi(B^2)_4 = \chi(O^2)_4 = -1$; consequently, the operators of the families $B^2$ and $O^2$ are not bijective. 

On the other hand, let $B \in B^1 \cup B^3 \cup O^1$; we show that $B$ and $\widetilde{B}$ are bijections. 
Let $B \in B^1 \cup O^1$; in this case, $k \neq 0, -1$ and $l_\alpha + m_\beta = m_\beta \neq -1$, namely, the operator $B$ is a bijection. 
By Lemma 5, the equality $\chi(\widetilde{B})_7 = -k - 1$ holds; therefore, 
$\chi(\widetilde{B})_7 \neq 0, -1$ since $k \neq 0, -1$. 
Finally, $\chi(\widetilde{B})_1 + \chi(\widetilde{B})_4 = -1 - m_\beta - 1$ by Lemma 5; 
hence, we obtain that $\chi(\widetilde{B})_1 + \chi(\widetilde{B})_4 \neq -1$ due to the fact that $m_\beta \neq -1$. Thus, the operator $\widetilde{B}$ is also bijective.

Now let $B \in B^3$, that is, $B = B^3_m$ for some $1 \leq m < (p-1)/2$; then $l_\alpha + m_\beta = 2 m \neq -1$ and $k = \frac{m^2}{2m+1}$. We find that $k \neq 0$ since $m \neq 0$. The equality $\frac{m^2}{2m+1} = -1$ is equivalent to $(m + 1)^2 = 0$, since $m \neq -1 = p-1$ by assumption, we obtain that $k \neq -1$. 
Consequently, $B$ is bijective. 
Similarly to the previous case, we obtain $\chi(\widetilde{B})_1 + \chi(\widetilde{B})_4 = -m-2$;
therefore, $\chi(\widetilde{B})_1 + \chi(\widetilde{B})_4 \neq -1$ since $m \neq -1$. 
Thus, 
$\chi(\widetilde{B})_7 = -k - 1 \neq 0, -1$ since $k = \frac{m^2}{2m+1} \neq 0, -1$.
As a result, we find that $\widetilde{B}$ is a bijection.
\hfill $\square$

\subsection{Construction for Bijective Operators}

\textbf{Proposition 12}.
Let $G$ be a group of nilpotency class $2$, i.\,e., $[[G, G], G] = \{ e \}$, and let $\varphi \in \Aut(G)$. Then the mapping $B(g) = \theta(g) \varphi(g)$, where $g \in G$ and $\theta \colon G \to Z(G)$, is a Rota---Baxter operator if and only if 
$\theta(g) \theta(h) = \theta(g h) \theta([\varphi(g^{-1}), h^{-1}])\varphi([\varphi(g^{-1}), h^{-1}])$.

{\sc Proof}. It is sufficient to expand by definition, assuming that $B$ is a Rota---Baxter operator:
\begin{multline*}
\theta(g) \theta(h) \varphi(g) \varphi(h) = B(g) B(h) = B(g B(g) h B(g)^{-1}) = B(g \varphi(g) h \varphi(g)^{-1}) \\
 = \theta(g \varphi(g) h \varphi(g)^{-1}) \varphi(g \varphi(g) h \varphi(g)^{-1} h^{-1} h) 
 = \theta(g \varphi(g) h \varphi(g)^{-1}) \varphi(g [\varphi(g)^{-1}, h^{-1}] h)  \\
 = \theta(g \varphi(g) h \varphi(g)^{-1}) \varphi(g) \varphi(h) \varphi([\varphi(g)^{-1}, h^{-1}]).
\end{multline*}
Thus, $B$ is a Rota---Baxter operator if and only if 
\begin{equation} \label{theta_eq}
\theta(g) \theta(h) = \theta(g \varphi(g) h \varphi(g)^{-1}) \varphi([\varphi(g)^{-1}, h^{-1}]).
\end{equation}

We choose $h = z \in Z(G)$ in formula \eqref{theta_eq} and obtain that $\theta(g)\theta(z) = \theta(gz)$. Then it is true that
\begin{multline*}
\theta(g) \theta(h) = \theta(g \varphi(g) h \varphi(g)^{-1}) \varphi([\varphi(g)^{-1}, h^{-1}]) 
 = \theta( g [\varphi(g)^{-1}, h^{-1}] h) \varphi([\varphi(g)^{-1}, h^{-1}]) \\
\qquad \qquad \qquad  \qquad = \theta(gh) \theta([\varphi(g)^{-1}, h^{-1}]) \varphi([\varphi(g)^{-1}, h^{-1}]) 
 = \theta(gh) B([\varphi(g)^{-1}, h^{-1}]). \hfill \square
\end{multline*}

\textbf{Proposition 13}. 
All bijective Rota---Baxter operators on the group $H_3$ have the form as in Proposition~11.

{\sc Proof}. 
Let $B$ be a Rota---Baxter operator on $H_3$ such that $\chi(B) = (l_\alpha, l_\beta, m_\alpha, m_\beta$, $x_\alpha, x_\beta, k)$ and let $B$ be a bijection. 
Then $l_\alpha m_\beta - l_\beta m_\alpha = k(1 + l_\alpha + m_\beta) \neq 0$ by \eqref{RB_condition} and the bijectivity condition. 
We consider $\varphi \in \Aut(G)$ with parameters $u_\alpha = l_\alpha$, $u_\beta = l_\beta$, $v_\alpha = m_\alpha$, $v_\beta = m_\beta$, $\Delta = u_\alpha v_\beta - u_\beta v_\alpha \neq 0$, and arbitrary $w_\alpha$, $w_\beta$. We set $\theta(b^s a^t c^r) = c^{B_\gamma(t, s, r) - \varphi(t, s, r)}$; then the equality $B(g) = \theta(g) \varphi(g)$ holds for any $g \in G$ by formulas \eqref{B_alpha_beta}, \eqref{B_gamma}, \eqref{phi_alpha_beta}, \eqref{phi_gamma}, and the operator $B$ has the form as in Proposition 12.
\hfill $\square$

\textbf{Example 2}.
Consider the Rota---Baxter operator $B(g) = g^{-1}$ on $H_3$. We choose an automorphism $\varphi(b^s a^t c^r) = b^{-s} a^{-t} c^r$ with parameters $u_\alpha = -1$, $u_\beta = 0$, $v_\alpha = 0$, $v_\beta = -1$, $w_\alpha = w_\beta = 0$, and 
$\theta(b^s a^t c^r) = c^{ts - 2r}$. 
We obtain $B(g) = \theta(g) \varphi(g)$ for any $g \in G$.

\textbf{Example 3}.
Consider the Rota---Baxter operator $B^1_{1, 1}(b^s a^t c^r) = b^{2t + s} a^{-s} c^{-ts + t(t-1) + r}$ on $H_3$. We choose an automorphism $\varphi(b^s a^t c^r) = b^{2t + s} a^{-s} c^{-2ts -\frac{s(s-1)}{2} + 2r}$ with parameters $u_\alpha = 0$, $u_\beta = -1$, $v_\alpha = 2$, $v_\beta = 1$, $w_\alpha = w_\beta = 0$, and we choose $\theta(b^s a^t c^r) = c^{ts + t(t-1) + \frac{s(s - 1)}{2} - r}$. We obtain that $B^1_{1, 1}(g) = \theta(g) \varphi(g)$ for any $g \in G$.

\textbf{Proposition 14}.
Let $G$ and $H$ be groups, let $\psi \colon Z(G) \to Z(H)$ be an isomorphism, and let $B_G$ be a Rota---Baxter operator on $G$ and $B_H$ be a Rota---Baxter operator on $H$ such that $B(Z(G)) \subseteq Z(G)$, $B(Z(H)) \subseteq Z(H)$, and $\psi(B_G(z)) = B_H(\psi(z))$ for any $z \in Z(G)$. 
Then the mapping $B \colon G \times_\psi H \to G \times_\psi H$ defined by the formula $B( (g, h) K ) = (B_G(g), B_H(h)) K$, where $K$ is defined as in \eqref{center_product}, is a Rota---Baxter operator on $G \times_\psi H$.

{\sc Proof}. We verify well-definedness. Since $B_G$ preserves the center, we have $B_G(zg) = B_G(z) B_G(g)$ for any $z \in Z(G)$, and a similar relation holds for $B_H$. Let $g \in G$, $h \in H$, and $z \in Z(G)$; then  
\begin{multline*}
B((g z, h \psi(z))K) 
 = (B_G(g) B_G(z), B_H(h) B_H(\psi(z))) K \\
 = (B_G(g), B_H(h)) (B_G(z), B_H(\psi(z))) K 
 = (B_G(g), B_H(h)) K = B((g, h)K).
\end{multline*}
Thus, the well-definedness is proved.
We show that $B$ is a Rota---Baxter operator:
\begin{multline*}
B((g_1, h_1)K) B((g_2, h_2)K) 
 = (B_G(g_1) B_G(g_2), B_H(h_1) B_H(h_2)) K  \\
 = (B_G(g_1 B_G(g_1) g_2 B_G(g_1)^{-1}), B_H(h_1 B_H(h_1) h_2 B_H(h_1)^{-1})) K \\
 = B((g_1, h_1) \cdot (B_G(g_1), B_H(h_1)) \cdot (g_2, h_2) \cdot (B_G(g_1), B_H(h_1))^{-1}K)  \\
 \qquad\qquad\qquad \qquad \qquad  \qquad = B((g_1, h_1) \cdot B((g_1, h_1)) \cdot (g_2, h_2) \cdot B((g_1, h_1))^{-1}K). 
 \hfill \square
\end{multline*}

A similar proposition holds for group automorphisms.

\textbf{Proposition 15}.
Let $G$ and $H$ be groups, let $\psi \colon Z(G) \to Z(H)$ be an isomorphism, and let $\varphi_G \in \Aut(G)$, $\varphi_H \in \Aut(H)$, and $\psi(\varphi_G(z)) = \varphi_H(\psi(z)))$ for any $z \in Z(G)$. 
Then the mapping $\varphi \colon G \times_\psi H \to G \times_\psi H$ 
defined by the formula $\varphi( (g, h) K ) = (\varphi_G(g), \varphi_H(h)) K$, 
where $K$ is defined by~\eqref{center_product}, is an automorphism of $G \times_\psi H$.

The proof of Proposition~15 is similar to the proof of Proposition 14.

\textbf{Corollary 1}.
On an extraspecial group $G$ of exponent $p$,
there exist Rota---Baxter operators $B$ such that $B$ and $\widetilde{B}$ are bijective and are defined according to Proposition~11.

{\sc Proof}.
By assumption, $G \cong H_3^n/K$, where $K = \{ (z, \ldots, z) \in Z(H_3)^n \}$. 
On each copy of~$H_3$, a bijective operator $B = B^1_{1, 1}$ is defined, for which there exist $\varphi \in \Aut(H_3)$ and a mapping $\theta \colon H_3 \to Z(H_3)$ such that $B(g) = \theta(g) \varphi(g)$ for any $g \in H_3$. 
We can choose $\varphi$ and $\theta$ from Example 4. We choose mappings $D, \Phi \colon G \to G$ such that 
$$
D( (g_1, \ldots, g_n) K ) = (B(g_1), \ldots, B(g_n)) K, \quad
\Phi((g_1, \ldots, g_n) K) = (\varphi(g_1), \ldots, \varphi(g_n)) K. 
$$
Then, by Propositions 13 and 14, the mapping $D$ is a Rota---Baxter operator on $G$ and $\Phi \in \Aut(G)$. 
We define a mapping $\Theta \colon G \to Z(G)$ by the formula $\Theta((g_1, \ldots, g_n) K) = (\theta(g_1), \ldots, \theta(g_n)) K$.

We obtain that $D(\bar{g}) = \Theta(\bar{g}) \Phi(\bar{g})$, where $\bar{g} \in G$. Moreover, since $\widetilde{D}( (g_1, \ldots, g_n) K ) = (\widetilde{B}(g_1), \ldots, \widetilde{B}(g_n)) K$, the operator $\widetilde{B}$ is also bijective by Proposition~11, and, consequently, satisfies the construction described in Proposition~12 by Proposition 13. 
Then, repeating the arguments applied above for $D$, we obtain that $\widetilde{D}$ is given according to the construction described in Proposition~12.

\subsection{Examples}

We consider the matrix representation of the Heisenberg group $H_3$ and we set
\begin{equation*}
a = 
\begin{pmatrix} 
 1 & 1 & 0 \\
 0 & 1 & 0 \\
 0 & 0 & 1
\end{pmatrix}, \quad
b = 
\begin{pmatrix} 
 1 & 0 & 0 \\
 0 & 1 & 1 \\
 0 & 0 & 1
\end{pmatrix}, \quad
c = 
\begin{pmatrix} 
 1 & 0 & 1 \\
 0 & 1 & 0 \\
 0 & 0 & 1
\end{pmatrix}
\in M_3(\mathbb{Z}_p).
\end{equation*}
The relations $[a, b] = c$ and $a^p = b^p = c^p = e$ hold, and
$$
B \left(
\begin{pmatrix} 
 1 & t & r \\
 0 & 1 & s \\
 0 & 0 & 1
\end{pmatrix}
\right)
= B(b^s a^t c^r) = \begin{pmatrix} 
 1 & B_\alpha(t, s) & B_\gamma(t, s, r) \\
 0 & 1 & B_\beta(t, s) \\
 0 & 0 & 1
\end{pmatrix},$$
where $B_\alpha, B_\beta, B_\gamma$ are the mappings defined as in \eqref{B_alpha_beta} and \eqref{B_gamma}.

The classification obtained in Theorem 1 covers the examples from the papers \cite{Ita}, \cite{Goncharovplus}, and \cite{Zhu}. A.~Caranti and L.~Stefanello, in Section 6 of their paper \cite{Ita}, provide an example of the following Rota---Baxter operator on $H_3(\mathbb{Z}_p)$:
$$
B \left(
\begin{pmatrix} 
 1 & t & r \\
 0 & 1 & s \\
 0 & 0 & 1
\end{pmatrix}
\right)
=
\begin{pmatrix} 
 1 & dt & \frac{d^3}{1+2d} ts + \frac{d^2}{1+2d} r \\
 0 & 1 & ds \\
 0 & 0 & 1
\end{pmatrix},
$$
where $d \neq -\frac{1}{2}$ is a parameter.
We have $\chi(B) = (d, 0, 0, d, 0, 0, \frac{d^2}{1+2d})$. 
If $d = 0$, then $B$ is trivial;
if $1 \leq d < \frac{p-1}{2}$, then $B = B^3_d$;
if $\frac{p-1}{2} < d \leq p-2$, then $\widetilde{B} = B^3_{p - d - 1}$.

In Section 4 of the paper \cite{Zhu}, Y. Jiang, Y. Sheng, and C. Zhu define a Rota---Baxter operator on $H_3(\mathbb{R})$ of the form
$$
B \left(
\begin{pmatrix} 
 1 & t & r \\
 0 & 1 & s \\
 0 & 0 & 1
\end{pmatrix}
\right)
=
\begin{pmatrix} 
 1 & 0 & s \\
 0 & 1 & 0 \\
 0 & 0 & 1
\end{pmatrix}.
$$

We consider this Rota---Baxter operator as an operator on $H_3(\mathbb{Z}_p)$. By item a) of Proposition 9, there exists $\varphi \in \Aut(G)$ such that $B^{(\varphi)} = B^4$.

In 2024, M. Goncharov, P. Kolesnikov, Y. Sheng, and R. Tang, in their paper \cite{Goncharovplus}, provided an example of a Rota---Baxter operator over $H_3(\mathbb{C})$:
$$
B \left(
\begin{pmatrix} 
 1 & t & r \\
 0 & 1 & s \\
 0 & 0 & 1
\end{pmatrix}
\right)
=
\begin{pmatrix} 
 1 & s & -\frac{t^2}2+\frac{s^2}2+ts-r \\
 0 & 1 & t \\
 0 & 0 & 1
\end{pmatrix}.
$$

We consider this Rota---Baxter operator as a Rota---Baxter operator on $H_3(\mathbb{Z}_p)$, $\chi(B) = (0, 1, 1, 0, -\frac{1}{2}, \frac{1}{2}, -1)$. Then 
$\chi(\widetilde{B}) = (-1, -1, -1, -1, \frac{1}{2}, \frac{1}{2}, 0)$ by Lemma 5. 
Then, choosing the automorphism $\varphi$ with parameters $t_\alpha = 2$, $t_\beta = 1$, $s_\alpha = -2$, and $s_\beta = 1$, 
we obtain by formula~\eqref{coeff_aut_H3} that
$\chi(\widetilde{B}^{(\varphi)}) = (0, 0, 0, -2, 0, x_\beta^{(\varphi)}, 0)$. 
By item c) of Proposition~9, there exists $\psi \in \Aut(G)$ such that $\widetilde{B}^{(\psi \varphi)} = B^8_{p-2}$.

\subsection{Splitting Operators on $H_3$ and $E_3$}

\textbf{Proposition 16}.

a) Let $B$ be a non-trivial splitting operator on the group $H_3$; then $B$, up to conjugation by automorphisms and the action \ $\widetilde{}$ \, coincides with the operator $B^8_{-1}$.

b) Let $B$ be a non-trivial splitting operator on the group $E_3$; then $B$, up to conjugation by automorphisms and the action \ $\widetilde{}$ \, coincides with the operator $\mathfrak{B}^4_{-1, 0}$.

{\sc Proof}.
We use Proposition 3 and Lemma 5 in all items of the proof.

a) If the operator is splitting, then $k = 0$ or $k = -1$; we consider the operators of the families $B^4, B^5, B^6, B^7, B^8$.

a1) $\chi(B^4) = (0, 0, 0, 0, 1, 0, 0)$, $\chi(\widetilde{B^4}) = (-1, 0, 0, -1, -1, 0, -1)$. 
In this case, $(B^4 \circ \widetilde{B^4})(a) = c^{-1}$.

a2) $\chi(B^5) = (0, 1, 0, 0, 1, 0, 0)$, $\chi(\widetilde{B^5}) = (-1, -1, 0, -1, -1, 0, -1)$. 
We obtain that $(B^5 \circ \widetilde{B^5})(a) = c^{-1}$.

a3) $\chi(B^6) = (0, 1, 0, 0, 0, 0, 0)$, $\chi(\widetilde{B^6}) = (-1, -1, 0, -1, 0, 0, -1)$. 
Hence $(B^6 \circ \widetilde{B^6})(b) = a^{-1}$.

a4) We choose some $m \in \mathbb{Z}_p^*$, $\chi(B^7_m) = (0, 0, 0, m, 1, 0, 0)$, $\chi(\widetilde{B^7_m}) = (-1, 0, 0, -(m + 1), -1, 0, -1)$. Then $(B^7_m \circ \widetilde{B^7_m})(a) = c^{-1}$.

a5) Consider the operator $B^8_{-1}$; then $\chi(B^8_{-1}) = (0, 0, 0, -1, 0, 0, 0)$ and $\chi(\widetilde{B^8_{-1}}) = (-1, 0, 0, 0,0, -1, 0)$. One can easily understand that $B^8_{-1}$ is splitting. In this case, we have a exact factorization $H_3 = \Imm(B) \Imm(\widetilde{B}) = \langle b \rangle \langle a, c \rangle$.
We choose some $m \in \mathbb{Z}_p^* \setminus \{ p-1 \}$; then $\chi(B^8_m) = (0, 0, 0, m, 0, 0, 0)$, $\chi(\widetilde{B^8_m}) = (-1, 0, 0, -(m + 1), 0, 0, -1)$. Consequently, $(B^8_m \circ \widetilde{B^8_m})(b) = b^{-m(m+1)} \neq e$.

b) We prove that all Rota---Baxter operators on $E_3$, except for the trivial one and $\mathfrak{B}^4_{-1, 0}$, are non-splitting.

b1) $\chi_E(\mathfrak{B}^1) = [0, 0, 0, 0, 1, 0, 0]$, $\chi_E(\widetilde{\mathfrak{B}^1}) = [-1, 0, 0, -1, -1, 0, -1]$. We obtain that $(\mathfrak{B}^1 \circ \widetilde{\mathfrak{B}^1})(a) = c^{-1}$.

b2) $\chi_E(\mathfrak{B}^2_x) = [0, 0, 0, 0, x, x, 0]$, $\chi_E(\widetilde{\mathfrak{B}^2_x}) = [-1, 0, 0, -1, -x, -x, -1]$ for $1 \leq x \leq p-1$. Then $(\mathfrak{B}^2_{x} \circ \widetilde{\mathfrak{B}^2_{x}})(a) = c^{-x}$.

b3) $\chi_E(\mathfrak{B}^3_x) = [1, 1, -1, -1, x, 0, 0]$, $\chi_E(\widetilde{\mathfrak{B}^3_x}) = [-2, -1, 1, 0, -(2+x), -1, -1]$ where $0 \leq x \leq p-1$. We obtain that $(\mathfrak{B}^3_{x} \circ \widetilde{\mathfrak{B}^3_{x}})(b) = ba^{-1}c^{-x} \neq e$.

b4) $\chi_E(\mathfrak{B}^4_{l, y}) = [l, 0, -l, 0, 0, y, 0]$, $\chi_E(\widetilde{\mathfrak{B}^4_{l, y}}) = [-(l+1), 0, l, -1, -l(l+1), -y, -1]$ where $1 \leq l \leq p-1$ and $0 \leq y \leq p-1$. 
Let $l \neq -1$ or $y \neq 0$; then 
$(\mathfrak{B}^4_{l, y} \circ \widetilde{\mathfrak{B}^4_{l, y}})(a) = b^{l(l+1)} a^{-l(l+1)} c^{y l - \frac{1}{2} l^2 (l + l)(l + 2)} \neq e$. 
For $l = -1$ and $y = 0$, we have $\chi_E(\mathfrak{B}^4_{-1, 0}) = [-1, 0, 1, 0, 0, 0, 0]$, $\chi_E(\widetilde{\mathfrak{B}^4_{-1, 0}}) = [0, 0, -1, -1, 0, 0, -1]$. We compute
$$
\mathfrak{B}^4_{-1, 0} \widetilde{\mathfrak{B}^4_{-1, 0}} (b^s a^t c^r) = \mathfrak{B}^4_{-1, 0} \left(b^{-t-s} c^{ts + \frac{t(t-1)}{2} - r} \right) = e.
$$
We find that the operator $\mathfrak{B}^4_{-1, 0}$ is splitting. 
In this case, we have an exact factorization $E_3 = \Imm(B) \Imm(\widetilde{B}) = \langle b a^{-1} \rangle \langle b, c \rangle$.
\hfill $\square$

\section{Acknowledgements}
The author is deeply grateful to his supervisor, Vsevolod Gubarev, for his support, the scientific guidance and assistance in preparing this paper.

\noindent Andrey Savelyev \\
Novosibirsk State University \\
Pirogova str. 1, 630090 Novosibirsk, Russia \\
e-mail: andrsav2zz2@gmail.com

\newpage

\appendix
\renewcommand{\thesection}{\Alph{section}}
\section{Calculation of $\varphi^{-1}$}
We consider $\varphi \in \Aut(G)$ with parameters $u_\alpha, u_\beta, v_\alpha, v_\beta, w_\alpha, w_\beta$ and automorphisms $\psi$  with parameters $u_\alpha', u_\beta', v_\alpha', v_\beta', w_\alpha', w_\beta'$, defined by formula \eqref{op_automorphism}. Then
\begin{gather*}
(\psi \varphi)_{\alpha}(t, s) = u_\alpha'(u_\alpha t + u_\beta s) + u_\beta'(v_\alpha t + v_\beta s) = \frac{1}{\Delta}( u_\alpha v_\beta t + u_\beta v_\beta s - u_\beta v_\alpha t - u_\beta v_\beta s ) = t, \\
(\psi \varphi)_{\beta}(t, s) = v_\alpha'(u_\alpha t + u_\beta s) + v_\beta'(v_\alpha t + v_\beta s) = \frac{1}{\Delta}( -u_\alpha v_\alpha t - u_\beta v_\alpha s + u_\alpha v_\alpha t + u_\alpha v_\beta s ) = s,
\end{gather*}

\begin{multline*}
\Delta^2((\psi \varphi)_{\gamma}(t, s, r) - r) = \Delta^2 \Big( w_\alpha'(u_\alpha t + u_\beta s) + w_\beta'(v_\alpha t + v_\beta s) + u_\beta' v_\alpha' (u_\alpha t + u_\beta s)(v_\alpha t + v_\beta s) \\
+ \frac{u_\alpha' v_\alpha'}{2} (u_\alpha t + u_\beta s) (u_\alpha t + u_\beta s - 1) + \frac{u_\beta' v_\beta'}{2}(v_\alpha t + v_\beta s) (v_\alpha t + v_\beta s - 1) \\
+ (u_\alpha' v_\beta' - u_\beta' v_\alpha') \big( w_\alpha t + w_\beta s + u_\beta v_\alpha ts + u_\alpha v_\alpha \frac{t(t-1)}{2} + u_\beta v_\beta \frac{s(s-1)}{2} + (u_\alpha v_\beta - u_\beta v_\alpha) r \big) \Big) - \Delta^2 r \\
= \Big(-w_\alpha v_\beta + w_\beta v_\alpha + \frac{u_\alpha v_\alpha v_\beta}{2} - \frac{u_\beta v_\alpha v_\beta}{2} - \frac{v_\alpha v_\beta}{2} \Big) (u_\alpha t + u_\beta s) \\
+ \Big( w_\alpha u_\beta - w_\beta u_\alpha -\frac{u_\alpha u_\beta v_\alpha}{2} + \frac{u_\alpha u_\beta v_\beta}{2} - \frac{u_\alpha u_\beta}{2} \Big)(v_\alpha t + v_\beta s) \\
+ u_\beta v_\alpha (u_\alpha t + u_\beta s)(v_\alpha t + v_\beta s) - \frac{v_\alpha v_\beta}{2}(u_\alpha t + u_\beta s)(u_\alpha t + u_\beta s - 1) - \frac{u_\alpha u_\beta}{2}(v_\alpha t + v_\beta s)(v_\alpha t + v_\beta s - 1) \\
+ \Delta \Big( w_\alpha t + w_\beta s + u_\beta v_\alpha ts + u_\alpha v_\alpha \frac{t(t-1)}{2} + u_\beta v_\beta \frac{s(s-1)}{2} \Big) \\
= -w_\alpha u_\alpha v_\beta t - w_\alpha u_\beta v_\beta s + w_\beta u_\alpha v_\alpha t + w_\beta u_\beta v_\alpha s + \frac{u_\alpha^2 v_\alpha v_\beta}{2} t + \frac{u_\alpha u_\beta v_\alpha v_\beta}{2} s - \frac{u_\alpha u_\beta v_\alpha v_\beta}{2} t - \frac{u_\beta^2 v_\alpha v_\beta}{2} s \\
- \frac{u_\alpha v_\alpha v_\beta}{2} t - \frac{u_\beta v_\alpha v_\beta}{2} s + w_\alpha u_\beta v_\alpha t + w_\alpha u_\beta v_\beta s - w_\beta u_\alpha v_\alpha t - w_\beta u_\alpha v_\beta s - \frac{u_\alpha u_\beta v_\alpha^2}{2} t - \frac{u_\alpha u_\beta v_\alpha v_\beta}{2} s \\
+ \frac{u_\alpha u_\beta v_\alpha v_\beta}{2} t + \frac{u_\alpha u_\beta v_\beta^2}{2} s -\frac{u_\alpha u_\beta v_\alpha}{2} t - \frac{u_\alpha u_\beta v_\beta}{2} s + u_\alpha u_\beta v_\alpha^2 t^2 + u_\alpha u_\beta v_\alpha v_\beta ts + u_\beta^2 v_\alpha^2 ts + u_\beta^2 v_\alpha v_\beta s^2 \\
- \frac{u_\alpha^2 v_\alpha v_\beta}{2} t^2 - u_\alpha u_\beta v_\alpha v_\beta ts + \frac{u_\alpha v_\alpha v_\beta}{2} t 
- \frac{u_\beta^2 v_\alpha v_\beta}{2} s^2 + \frac{u_\beta v_\alpha v_\beta}{2} s - \frac{u_\alpha u_\beta v_\alpha^2}{2} t^2 - u_\alpha u_\beta v_\alpha v_\beta ts \\
+ \frac{u_\alpha u_\beta v_\alpha}{2} t
- \frac{u_\alpha u_\beta v_\beta^2}{2} s^2 + \frac{u_\alpha u_\beta v_\beta}{2} s + w_\alpha u_\alpha v_\beta t - w_\alpha u_\beta v_\alpha t + w_\beta u_\alpha v_\beta s - w_\beta u_\beta v_\alpha s \\
+ u_\alpha u_\beta v_\alpha v_\beta ts - u_\beta^2 v_\alpha^2 ts + \frac{u_\alpha^2 v_\alpha v_\beta}{2} t^2 - \frac{u_\alpha u_\beta v_\alpha^2}{2} t^2 - \frac{u_\alpha^2 v_\alpha v_\beta}{2} t + \frac{u_\alpha u_\beta v_\alpha^2}{2} t + \frac{u_\alpha u_\beta v_\beta^2}{2} s^2 - \frac{u_\beta^2 v_\alpha v_\beta}{2} s^2 \\
- \frac{u_\alpha u_\beta v_\beta^2}{2} s + \frac{u_\beta^2 v_\alpha v_\beta}{2} s
= \Big( u_\alpha u_\beta v_\alpha v_\beta + u_\beta^2 v_\alpha^2 
- u_\alpha u_\beta v_\alpha v_\beta - u_\alpha u_\beta v_\alpha v_\beta + u_\alpha u_\beta v_\alpha v_\beta - u_\beta^2 v_\alpha^2 \Big) ts \\
+ \Big( u_\alpha u_\beta v_\alpha^2 - \frac{u_\alpha^2 v_\alpha v_\beta}{2} - \frac{u_\alpha u_\beta v_\alpha^2}{2} + \frac{u_\alpha^2 v_\alpha v_\beta}{2} - \frac{u_\alpha u_\beta v_\alpha^2}{2} \Big) t^2 \\
+ \Big( u_\beta^2 v_\alpha v_\beta - \frac{u_\beta^2 v_\alpha v_\beta}{2} - \frac{u_\alpha u_\beta v_\beta^2}{2} + \frac{u_\alpha u_\beta v_\beta^2}{2} - \frac{u_\beta^2 v_\alpha v_\beta}{2} \Big) s^2 \\
+ \Big( -w_\alpha u_\alpha v_\beta + w_\beta u_\alpha v_\alpha + \frac{u_\alpha^2 v_\alpha v_\beta}{2} - \frac{u_\alpha u_\beta v_\alpha v_\beta}{2} - \frac{u_\alpha v_\alpha v_\beta}{2} + w_\alpha u_\beta v_\alpha - w_\beta u_\alpha v_\alpha \allowdisplaybreaks \\
- \frac{u_\alpha u_\beta v_\alpha^2}{2} + \frac{u_\alpha u_\beta v_\alpha v_\beta}{2} - \frac{u_\alpha u_\beta v_\alpha}{2} + \frac{u_\alpha v_\alpha v_\beta}{2} + \frac{u_\alpha u_\beta v_\alpha}{2} + w_\alpha u_\alpha v_\beta - w_\alpha u_\beta v_\alpha - \frac{u_\alpha^2 v_\alpha v_\beta}{2} + \frac{u_\alpha u_\beta v_\alpha^2}{2} \Big) t \\
+ \Big( -w_\alpha u_\beta v_\beta + w_\beta u_\beta v_\alpha + \frac{u_\alpha u_\beta v_\alpha v_\beta}{2} - \frac{u_\beta^2 v_\alpha v_\beta}{2} - \frac{u_\beta v_\alpha v_\beta}{2} + w_\alpha u_\beta v_\beta - w_\beta u_\alpha v_\beta - \frac{u_\alpha u_\beta v_\alpha v_\beta}{2} \\
+ \frac{u_\alpha u_\beta v_\beta^2}{2}
- \frac{u_\alpha u_\beta v_\beta}{2} + \frac{u_\beta v_\alpha v_\beta}{2} + \frac{u_\alpha u_\beta v_\beta}{2} + w_\beta u_\alpha v_\beta - w_\beta u_\beta v_\alpha - \frac{u_\alpha u_\beta v_\beta^2}{2} + \frac{u_\beta^2 v_\alpha v_\beta}{2} \Big) s = 0.
\end{multline*}
We showed that $\psi = \varphi^{-1}$.

\section{Calculation $B^{(\varphi)}$}
Let $B$ be RB-operator with parameters $l_\alpha, l_\beta, m_\alpha, m_\beta, x_\alpha, x_\beta, k$.
Let us consider $B^{(\varphi)}$ and find the set of parameters $(l_\alpha^{(\varphi)}, l_\beta^{(\varphi)}, m_\alpha^{(\varphi)}, m_\beta^{(\varphi)}, x_\alpha^{(\varphi)}, x_\beta^{(\varphi)}, k^{(\varphi)}) = \chi (B^{(\varphi)})$ for the automorphism $\varphi$ with parameters $u_\alpha, u_\beta, v_\alpha, v_\beta,w_\alpha,w_\beta$. We calculate
\begin{multline*}
(B^{(\varphi)})_\alpha(t, s) = u_\alpha'(l_\alpha (u_\alpha t + u_\beta s) + l_\beta(v_\alpha t + v_\beta s)) +  u_\beta'(m_\alpha (u_\alpha t + u_\beta s) + m_\beta (v_\alpha t + v_\beta s)) \\
= \frac{1}{\Delta} (u_\alpha v_\beta l_\alpha t + u_\beta v_\beta l_\alpha s + v_\alpha v_\beta l_\beta t + v_\beta^2 l_\beta s -  u_\alpha u_\beta m_\alpha t - u_\beta^2 m_\alpha s - u_\beta v_\alpha m_\beta t - u_\beta v_\beta m_\beta s) \\
= \frac{1}{\Delta}(u_\alpha v_\beta l_\alpha + v_\alpha v_\beta l_{\beta} - u_\alpha u_\beta m_\alpha - u_\beta v_\alpha m_\beta) t + \frac{1}{\Delta} (u_\beta v_\beta l_\alpha + v_\beta^2 l_\beta - u_\beta^2 m_\alpha - u_\beta v_\beta m_\beta) s,
\end{multline*}
\begin{multline*}
(B^{(\varphi)})_\beta(t, s) = v_\alpha'(l_\alpha (u_\alpha t + u_\beta s) + l_\beta(v_\alpha t + v_\beta s)) +  v_\beta'(m_\alpha (u_\alpha t + u_\beta s) + m_\beta(v_\alpha t + v_\beta s)) \\
= \frac{1}{\Delta}(-u_\alpha v_\alpha l_\alpha t - u_\beta v_\alpha l_\alpha s - v_\alpha^2 l_\beta t - v_\alpha v_\beta l_\beta s + u_\alpha^2 m_\alpha t + u_\alpha u_\beta m_\alpha s + u_\alpha v_\alpha m_\beta t + u_\alpha v_\beta m_\beta s) \\
= \frac{1}{\Delta}(-u_\alpha v_\alpha l_\alpha - v_\alpha^2 l_\beta + u_\alpha^2 m_\alpha + u_\alpha v_\alpha m_\beta) t + \frac{1}{\Delta}(-u_\beta v_\alpha l_\alpha - v_\alpha v_\beta l_\beta + u_\alpha u_\beta m_\alpha + u_\alpha v_\beta m_\beta) s.
\end{multline*}

Then $l_\alpha^{(\varphi)}, l_\beta^{(\varphi)}, m_\alpha^{(\varphi)}, m_\beta^{(\varphi)}$ have the form as in formulas \eqref{coeff_aut_H3}.
Next, we expand
\begin{multline*}
(B^{(\varphi)})_\gamma(t, s, r) = 
w_\alpha' (l_\alpha (u_\alpha t + u_\beta s) + l_\beta(v_\alpha t + v_\beta s)) + w_\beta' (m_\alpha (u_\alpha t + u_\beta s) + m_\beta(v_\alpha t + v_\beta s)) \\
+ u_\beta' v_\alpha' (l_\alpha (u_\alpha t + u_\beta s) + l_\beta(v_\alpha t + v_\beta s)) (m_\alpha (u_\alpha t + u_\beta s) + m_\beta(v_\alpha t + v_\beta s)) \\
+ \frac{1}{2} u_\alpha' v_\alpha' (l_\alpha (u_\alpha t + u_\beta s) + l_\beta(v_\alpha t + v_\beta s))(l_\alpha (u_\alpha t + u_\beta s) + l_\beta(v_\alpha t + v_\beta s) - 1) \\
+ \frac{1}{2} u_\beta' v_\beta' (m_\alpha (u_\alpha t + u_\beta s) + m_\beta(v_\alpha t + v_\beta s))(m_\alpha (u_\alpha t + u_\beta s) + m_\beta(v_\alpha t + v_\beta s) - 1) \\
+ (u_\alpha' v_\beta' - u_\beta' v_\alpha') \bigg((l_\beta m_\alpha + k m_\beta) (u_\alpha t + u_\beta s) (v_\alpha t + v_\beta s) + \frac{1}{2} (u_\alpha t + u_\beta s) (u_\alpha t + u_\beta s - 1) m_\alpha (l_\alpha + k) \\
+ \frac{1}{2} (v_\alpha t + v_\beta s) (v_\alpha t + v_\beta s - 1) l_\beta (m_\beta - k)
+ x_\alpha (u_\alpha t + u_\beta s) + x_\beta (v_\alpha t + v_\beta s) \\
+ k (w_\alpha t + w_\beta s + u_\beta v_\alpha ts + \frac{t(t-1)}{2} u_\alpha v_\alpha + \frac{s(s-1)}{2} u_\beta v_\beta + (u_\alpha v_\beta - u_\beta v_\alpha) r)\bigg) \allowdisplaybreaks \\
=
(l_\beta^{(\varphi)} m_\alpha^{(\varphi)} + k m_\beta^{(\varphi)} ) ts + m_\alpha^{(\varphi)} (l_\alpha^{(\varphi)} + k) \frac{t(t-1)}{2} + l_\beta^{(\varphi)} (m_\beta^{(\varphi)} - k) \frac{s(s-1)}{2} + R,
\end{multline*}
where $R$ is defined as
\begin{multline*}
R = (B^{(\varphi)})_\gamma(t, s, r)
 - (l_\beta^{(\varphi)} m_\alpha^{(\varphi)} + k m_\beta^{(\varphi)} ) ts 
 - m_\alpha^{(\varphi)} (l_\alpha^{(\varphi)} + k) \frac{t(t-1)}{2} 
 - l_\beta^{(\varphi)} (m_\beta^{(\varphi)} - k) \frac{s(s-1)}{2} \\
 =: R_{\alpha \beta} ts + R_{\alpha^2} t^2 + R_{\beta^2} s^2 + R_\alpha t + R_\beta s + R_\gamma r.
\end{multline*}

We have 
$R_\gamma = (u_\alpha v_\beta - u_\beta v_\alpha)(u_\alpha' v_\beta' - u_\beta' v_\alpha')k = k$, 
therefore, $k^{(\varphi)} = k$. 
Due to \eqref{B_gamma} and definition $R$, we obtain
$R_{\alpha \beta} = R_{\alpha^2} = R_{\beta^2} = 0$ and 
$x_\alpha^{(\varphi)} = R_\alpha$, $x_\beta^{(\varphi)} = R_\beta$.

Now we write the formulas for $R_\alpha$ and $R_\beta$:
\begin{multline}
R_\alpha = \frac{1}{2} m_\alpha^{(\varphi)} (l_\alpha^{(\varphi)} + k) + w_\alpha'(u_\alpha l_\alpha + v_\alpha l_\beta) + w_\beta'(u_\alpha m_\alpha + v_\alpha m_\beta) \\
- \frac{1}{2} u_\alpha' v_\alpha' (l_\alpha u_\alpha + l_\beta v_\alpha) - \frac{1}{2} u_\beta' v_\beta' (m_\alpha u_\alpha + m_\beta v_\alpha) \\
+ (u_\alpha' v_\beta' - u_\beta' v_\alpha') (- \frac{1}{2} u_\alpha m_\alpha (l_\alpha + k) - \frac{1}{2} v_\alpha l_\beta (m_\beta - k) + u_\alpha x_\alpha + v_\alpha x_\beta + w_\alpha k - \frac{1}{2} u_\alpha v_\alpha k)  \allowdisplaybreaks \\
= \frac{1}{2 \Delta^2}(-u_\alpha v_\alpha l_\alpha - v_\alpha^2 l_\beta + u_\alpha^2 m_\alpha + u_\alpha v_\alpha m_\beta)(u_\alpha v_\beta l_\alpha + v_\alpha v_\beta l_\beta - u_\alpha u_\beta m_\alpha - u_\beta v_\alpha m_\beta + \Delta \cdot k) \\
+ \frac{1}{2 \Delta^2} (-2 w_\alpha v_\beta + 2 w_\beta v_\alpha + u_\alpha v_\alpha v_\beta - u_\beta v_\alpha v_\beta ) (u_\alpha l_\alpha + v_\alpha l_\beta) \\
+ \frac{1}{2 \Delta^2} (2 w_\alpha u_\beta - 2 w_\beta u_\alpha - u_\alpha u_\beta v_\alpha + u_\alpha u_\beta v_\beta - ) (u_\alpha m_\alpha + v_\alpha m_\beta) \\
+ \frac{1}{2 \Delta^2} \cdot \Delta (-u_\alpha m_\alpha (l_\alpha + k) - v_\alpha l_\beta (m_\beta - k) + 2 u_\alpha x_\alpha + 2 v_\alpha x_\beta + 2 w_\alpha k - u_\alpha v_\alpha k); \label{R_a}
\end{multline}

\vspace{-1.0cm}

\begin{multline}
R_\beta = \frac{1}{2} l_\beta^{(\varphi)} (m_\beta^{(\varphi)} - k) + w_\alpha'(u_\beta l_\alpha + v_\beta l_\beta) + w_\beta'(u_\beta m_\alpha + v_\beta m_\beta) \\
- \frac{1}{2} u_\alpha' v_\alpha' ( u_\beta l_\alpha + v_\beta l_\beta ) - \frac{1}{2} u_\beta' v_\beta'( u_\beta m_\alpha + v_\beta m_\beta) \\
+ (u_\alpha' v_\beta' - u_\beta' v_\alpha') \Big( -\frac{1}{2} u_\beta m_\alpha (l_\alpha + k) - \frac{1}{2} v_\beta l_\beta (m_\beta - k) + u_\beta x_\alpha + v_\beta x_\beta + w_\beta k - \frac{1}{2} u_\beta v_\beta k \Big) \\
= \frac{1}{2 \Delta^2} (u_\beta v_\beta l_\alpha + v_\beta^2 l_\beta - u_\beta^2 m_\alpha - u_\beta v_\beta m_\beta)(-u_\beta v_\alpha l_\alpha - v_\alpha v_\beta l_\beta + u_\alpha u_\beta m_\alpha + u_\alpha v_\beta m_\beta - \Delta \cdot k) \\
+ \frac{1}{2 \Delta^2}(-2 w_\alpha v_\beta + 2 w_\beta v_\alpha + u_\alpha v_\alpha v_\beta - u_\beta v_\alpha v_\beta ) (u_\beta l_\alpha + v_\beta l_\beta) \\
+ \frac{1}{2 \Delta^2}(2 w_\alpha u_\beta - 2 w_\beta u_\alpha - u_\alpha u_\beta v_\alpha + u_\alpha u_\beta v_\beta ) (u_\beta m_\alpha + v_\beta m_\beta) \\
+ \frac{1}{2 \Delta^2} \cdot \Delta (-u_\beta m_\alpha (l_\alpha + k) - v_\beta l_\beta (m_\beta - k) + 2 u_\beta x_\alpha + 2 v_\beta x_\beta + 2 w_\beta k - u_\beta v_\beta k). \label{R_b}
\end{multline}

\end{document}